\documentclass[preprint,12pt]{elsarticle}

\usepackage{amssymb}
\usepackage{amsmath} 
\usepackage{graphicx}
\usepackage{hyperref}
\usepackage{booktabs}
\usepackage{listings}
\usepackage{xcolor}
\usepackage{ragged2e}
\usepackage{float}
\usepackage{algorithm} 
\usepackage{algpseudocode} 
\usepackage{caption}
\usepackage{setspace}
\usepackage[utf8]{inputenc}
\usepackage{lineno} 
\usepackage{xcolor}
\journal{}

\begin{document}

\begin{frontmatter}



\title{Steady-state and transient thermal stress analysis using a polygonal finite element method}

\author[inst1]{Yang Yang}
\affiliation[inst1]{organization={PowerChina Kunming Engineering Corporation Limited},
            city={Kunming},
            postcode={650051}, 
            state={Yunnan},
            country={China}}

\author[inst2]{Mingjiao Yan}

\affiliation[inst2]{organization={College of Water Conservancy and Hydropower Engineering, Hohai University},
            city={Nanjing},
            postcode={210098}, 
            state={Jiangsu},
            country={China}}
\author[inst2]{Zongliang Zhang}
\author[inst2]{Dengmiao Hao}
\author[inst3]{Xuedong Chen}
\author[inst1]{Weixiong Chen}
\affiliation[inst3]{organization={Lijiang Water Resources and Hydropower Survey, Design and Research Institute Co., Ltd,},
            city={Lijiang},
            postcode={674100}, 
            state={Yunnan},
            country={China}}

\begin{abstract}
This work develops a polygonal finite element method (PFEM) for the analysis of steady-state and transient thermal stresses in two-dimensional continua. The method employs Wachspress rational basis functions to construct conforming interpolations over arbitrary convex polygonal meshes, providing enhanced geometric flexibility and accuracy in capturing complex boundary conditions and heterogeneous material behavior. A quadtree-based acceleration strategy is introduced to significantly reduce computational cost through the reuse of precomputed stiffness and mass matrices. The PFEM is implemented in ABAQUS via a user-defined element (UEL) framework. Comprehensive benchmark problems, including multi-scale and non-matching mesh scenarios, are conducted to verify the accuracy, convergence properties, and computational efficiency of the method. Results indicate that the proposed PFEM offers notable advantages over conventional FEM in terms of mesh adaptability, solution quality, and runtime performance. The method shows strong potential for large-scale simulations involving thermal–mechanical coupling, complex geometries, and multi-resolution modeling.
\end{abstract}


\begin{highlights}
\item A polygonal FEM (PFEM) using Wachspress functions is developed for thermal stress analysis.
\item The method supports arbitrary convex polygons and naturally handles hanging nodes.
\item A quadtree-based acceleration technique reuses parent matrices to reduce computation.
\item The proposed method is implemented in ABAQUS UEL and validated through benchmark tests.
\item The PFEM exhibits improved convergence and efficiency compared to classical FEM.
\end{highlights}


\end{frontmatter}


\section{Introduction}
Thermal stress analysis plays a vital role in engineering and materials science, particularly in the design and assessment of structures subjected to temperature variations~\cite{fang2021thermal,li2022accurate,yin2024understanding,lu2023numerical}. As materials expand or contract under thermal loads, stress concentrations may develop, potentially resulting in structural failure if not properly managed. Accurate modeling and prediction of thermal stresses are therefore essential to ensure the reliability and safety of engineering systems~\cite{moreira2022accurate}. Although the theory of thermoelasticity is well established and analytical solutions are available for a limited set of idealized problems, practical engineering structures often involve complex geometries, heterogeneous materials, and anisotropic behavior, rendering analytical approaches largely impractical~\cite{wu2010edge,Shibahara2011}.

To overcome these challenges, numerical methods have become indispensable tools for analyzing thermoelastic problems. Common approaches include the finite element method (FEM)~\cite{luo2018survey,giunta2016thermal}, finite difference method (FDM)~\cite{chen2013numerical}, finite cell method (FCM)~\cite{Zander2012}, and boundary element method (BEM)~\cite{cheng2016thermal}. Among these, FEM remains the most widely used due to its robustness and versatility, serving as the backbone of thermal stress analysis in both academic research and industrial practice~\cite{nemat2009elastic,faizin2023structural}. For instance, Najibi et al.~\cite{najibi2021transient} examined nonlinear transient thermal stress behavior in a thick hollow functionally graded cylinder, while Wang et al.~\cite{wang2023analysis} investigated thermal stress distributions in coatings under static and dynamic thermal loading. These studies collectively underscore the efficacy of FEM in capturing complex thermal stress phenomena with high fidelity.

Despite its broad adoption, conventional FEM exhibits several well-known limitations. The accuracy and convergence behavior of FEM solutions are strongly influenced by mesh quality, particularly in the presence of geometric complexities or evolving domains~\cite{kumbhar2020development,li2020n}. In two-dimensional settings, triangular and quadrilateral elements are predominantly employed. While triangular elements offer geometric flexibility, they often compromise accuracy; conversely, quadrilateral elements typically provide higher accuracy but are less adaptable to complex boundaries, complicating mesh generation~\cite{song2009scaled,ye2021free}. Although mesh refinement strategies can improve accuracy in regions of interest~\cite{huynh2020polytree,jansari2019adaptive}, they may introduce distortions or discontinuities at the interfaces between coarse and refined regions, thereby degrading solution quality.

To address these limitations, a variety of alternative numerical frameworks have been developed, including meshfree methods~\cite{chen2017meshfree}, physics-informed neural networks (PINNs)~\cite{cai2021physics,raissi2019physics,luo2023novel}, and isogeometric analysis (IGA)~\cite{cottrell2009isogeometric,willems2024isogeometric}. In addition, polygonal element-based methods have garnered increasing attention due to their enhanced geometric flexibility and approximation capabilities. Notable examples include the polygonal scaled boundary finite element method (PSBFEM)~\cite{chen2017novel,yang2022novel}, polygonal smoothed finite element method (PSFEM)~\cite{liu2021cell,yan2024fast}, polygonal virtual element method (PVEM)~\cite{benvenuti2022extended}, and polygonal finite element method (PFEM)~\cite{biabanaki2014polygonal,nguyen2017polytree,wu2023polygonal}. Unlike traditional elements with fixed topologies, polygonal elements can possess an arbitrary number of sides, naturally supporting a greater number of nodes in interpolation schemes~\cite{ooi2014scaled}, which enhances both accuracy and convergence. Among these approaches, PFEM has gained particular prominence. It retains the variational foundation of classical FEM while extending its applicability to general polygonal meshes. This enables seamless integration with existing FEM infrastructure and simplifies implementation. Moreover, PFEM offers a flexible and efficient platform for high-fidelity simulations involving complex or evolving geometries.

In large-scale simulations, mesh refinement is typically applied in localized regions to capture gradients more accurately~\cite{huynh2020polytree,jansari2019adaptive}. However, such localized refinement often introduces mesh incompatibilities at subdomain interfaces, commonly referred to as hanging nodes. To mitigate this, non-matching mesh techniques have been proposed, allowing for different resolutions across subdomains~\cite{bitencourt2015coupling,manzoli2018adaptive,li2020n}. While effective, these techniques often necessitate special treatments to enforce continuity across elements.

A key advantage of PFEM is its ability to address the hanging node issue in a straightforward manner. By treating hanging nodes as additional vertices within polygonal elements~\cite{rajagopal2018hyperelastic,tabarraei2008extended}, PFEM eliminates the need for constraint enforcement or complex remeshing procedures. This feature enables seamless handling of non-matching interfaces, particularly in quadtree or adaptive meshes, by naturally incorporating irregular node configurations into the interpolation framework of PFEM. This capability significantly enhances meshing flexibility, improves numerical stability, and facilitates accurate modeling of irregular or time-evolving domains. As a result, PFEM emerges as a robust and practical framework for solving thermo-mechanical problems in arbitrarily complex geometries.

The remainder of this paper is organized as follows. Section~\ref{sec:theory} presents the theoretical formulation of PFEM, including the treatment of hanging nodes and the construction of polygonal elements using Wachspress basis functions. Section~\ref{sec:implementation} details the numerical implementation and integration of PFEM into ABAQUS via a user-defined element (UEL). Section~\ref{sec:examples} validates the proposed method through numerical examples, highlighting its accuracy, efficiency, and convergence performance in thermal stress simulations. Finally, Section~\ref{sec:conclusions} concludes the paper and outlines potential future research directions for advancing PFEM in nonlinear and multi-physics applications.

\section{Thermal stress analysis of the PFEM}
In this section, we developed a PFEM formulation for steady-state and transient thermal problems. The method is based on Wachspress rational basis functions, which satisfy essential properties such as partition of unity, linear precision, and Kronecker delta property. Furthermore, a triangulation-based numerical integration scheme is adopted to facilitate element-level matrix computation. The resulting framework allows seamless incorporation of non-standard elements and supports hanging nodes without additional constraints. The governing equations for thermal conduction and thermally induced stress are derived and discretized in detail below.
\label{sec:theory}
\subsection{Thermal analysis}
The governing equation of transient heat conduction in a two-dimensional (2D) domain $\Omega$ can be written as follows \cite{li2016novel}:
\begin{equation}
    \label{eq:gov1}
    \rho c \dot{\phi}-\nabla(\mathbf{k} \nabla \phi)-Q=0,
\end{equation}
where $\phi$ is the temperature, $\dot{\phi}$ is the temperature change rate, $Q$ is the rate of internal heat generation, $\rho$ is the density, $c$ is the specific heat. The thermal conductivity, denoted as $k$, is expressed as
\begin{equation}
    k=\begin{bmatrix}k_x&0\\0&k_y\end{bmatrix},
\end{equation}
where $k_x$ and $k_y$ are the thermal conductivity in the $x$ and $y$ directions, respectively.

With the initial conditions
\begin{equation}
\phi(x,y,t=0) =\phi_0(x,y) \quad \text{in} \quad \Omega,  
\end{equation}
and the boundary conditions
\begin{eqnarray}
    \phi(x,y)=\bar{\phi} \quad \text{in} \quad S_1,
    \label{eq:b1}
\end{eqnarray}
\begin{equation}
    -k \frac{\partial \phi}{\partial n}=\bar{q}_n= q_2 \quad \text { on } \quad S_2,
        \label{eq:b2}
\end{equation}
\begin{equation}
    -k \frac{\partial \phi}{\partial n}=\bar{q}_n= g\left(\phi-\phi_{\infty}\right) \quad \text { on } \quad S_3,
        \label{eq:b3}
\end{equation}
where $n$ is the outward normal vector to domain $\Omega$. $S$ is the boundary. Furthermore, $\bar{\phi}$ and $q$ are the prescribed boundary temperature and heat flux, respectively. $g$ is the convection heat transfer coefficient and $\phi_{\infty}$ is the ambient temperature.

The weak form is obtained by multiplying Eq. (\ref{eq:gov1}) with a test function $\delta w$ and integrating by parts:
\begin{equation}
    \begin{aligned}
    &\int_{\Omega} \rho c \, \delta w^\mathrm{T} \frac{\partial \phi}{\partial t} \, d\Omega + \int_{\Omega} (\mathbf{\nabla} \delta w)^\mathrm{T} k (\mathbf{\nabla} \phi) \, d\Omega - \int_{\Omega} \delta w^\mathrm{T} Q \, d\Omega \\
    &+ \int_{\Gamma_q} \delta w^\mathrm{T} q \, d\Gamma + \int_{\Gamma_h} \delta w^\mathrm{T} h (\phi - \phi_{\infty}) \, d\Gamma = 0\label{eq:galerkin}
    \end{aligned},
\end{equation}
where $\nabla$ denotes a differential operator defined as follows:
\begin{equation}
    \mathbf{\nabla}=\begin{bmatrix}\frac{\partial}{\partial x}\\\frac{\partial}{\partial y}\end{bmatrix}.
\end{equation}

In the Galerkin weak form given in Eq. (\ref{eq:galerkin}), the temperature field $\phi$ can be approximated as:
\begin{equation}
    \phi=\sum_{i=1}^mN_i \phi_i,
    \label{eq:temperature field}
\end{equation}
where \( \phi_i \) represents the temperature at node \( i \), and \( N_i \) is the corresponding shape function. By substituting Eq. (\ref{eq:temperature field}) into Eq. (\ref{eq:galerkin}), Eq. (\ref{eq:galerkin}) can then be expressed as
\begin{equation}
    \begin{aligned}
    &\int_{\Omega} \rho c \, N^{T} N \frac{\partial \mathbf{\phi}}{\partial t} \, d\Omega + \int_{\Omega} k\left[\frac{\partial N}{\partial x} ^{T} \frac{\partial N}{\partial x} + \frac{\partial N}{\partial y} ^{T}\frac{\partial N}{\partial y} \right] \mathbf{\phi} \, d\Omega \\ 
    &+ \int_{\Gamma_2} N^T q \, d\Gamma + \int_{\Gamma_3} {gN^T N}\phi \, d\Gamma - g \int_{\Gamma_3} N^T T_a \, d\Gamma - \int_{\Omega} N^T Q \, d\Omega \\&= \underbrace{\int_{\Omega} \rho c \, N^{T} N  \, d\Omega }_{\mathrm{M}}\frac{\partial \mathbf{\phi}}{\partial t}+\underbrace{\int_{\Omega} k\left[\frac{\partial N}{\partial x} ^{T} \frac{\partial N}{\partial x} + \frac{\partial N}{\partial y} ^{T}\frac{\partial N}{\partial y} \right]d\Omega}_{\mathrm{K}_{\mathrm{th}}}\phi \\&
    +\underbrace{\int_{\Gamma_2} N^T q \, d\Gamma-g \int_{\Gamma_3} N^T T_a \, d\Gamma - \int_{\Omega} N^T Q \, d\Omega}_{Q}=0.
    \end{aligned}
\end{equation}

The equilibrium equation of the discretized system can ultimately be written in the following matrix form:
\begin{equation}
\mathbf{M}\dot{\boldsymbol{\phi}}+\mathbf{K}_{\mathrm{th}}\boldsymbol{\phi}=\mathbf{Q},\label{eq:govingeq2}
\end{equation}
where $\mathbf{K}_{\mathrm{th}}$ is the thermal stiffness matrix. 

In which 
\begin{equation}
\mathbf{K}_{\mathrm{th}}=k\int_\Omega\mathbf{B}_{\mathrm{th}}^T\mathbf{B}_{\mathrm{th}}d\Omega,
\end{equation}
where $\mathbf{B}_{\mathrm{th}}$, referred to as the "strain-thermal" matrix, is given by  
\begin{equation}
    \mathbf{B}_\mathrm{th}=\nabla N=\begin{bmatrix}N_{i,x}\\N_{i,y}\end{bmatrix}.
\end{equation}

$\mathbf{M}_\mathrm{th}$ in Eq. (\ref{eq:govingeq2}) denotes thermal mass matrix, which can be written as
\begin{equation}
    \mathbf{M}_\mathrm{th}=\rho c\int_\Omega N^\mathrm{T} Nd\Omega.
\end{equation}

In this work, the backward difference method \cite{zienkiewicz2000finite} is employed to solve Eq.~(\ref{eq:govingeq2}). The time domain is discretized into several time steps, and the solution at each time node is progressively computed starting from the initial conditions. The nodal temperature at any intermediate time is then obtained through interpolation.

At time $[t,t+\Delta t]$, the temperature change rate $\dot{\boldsymbol{\phi}}$ can be expressed as 
\begin{equation}
    \dot{\boldsymbol{\phi}}(t)=\frac{\Delta \boldsymbol{\phi}}{\Delta t}=\frac{\boldsymbol{\phi}(t)^{t+\Delta t}-\boldsymbol{\phi}(t)^t}{\Delta t}. \label{eq:temperaturechange}
\end{equation}

Eq. (\ref{eq:temperaturechange}) is substituted into Eq. (\ref{eq:govingeq2}), and the equation at time step $t+\Delta t$ can be obtained as follows:
\begin{equation}
    (\mathbf{K}_\mathrm{th}+\frac{\mathbf{M}_\mathrm{th}}{\Delta t}) \boldsymbol{\phi}(t) =\mathbf{Q}+(\frac{\mathbf{M}_\mathrm{th}}{\Delta t}) \boldsymbol{\phi}(t-\Delta t).
\end{equation}

\subsection{Thermal mechanical analysis}
The governing equation describing linear elasto-static behavior subjected to thermal load is expressed as
\begin{equation}
    \label{eq:stressgov}
    \mathbf{\sigma}=\mathbf{D}(\mathbf{\varepsilon}-\mathbf{\varepsilon_0}),
\end{equation}
where $\mathbf{D}$ is the elastic matrix, and $\mathbf{\varepsilon_0}$ represents the initial strain caused by the temperature variation $\phi$. For 2D problems, the initial strain $\mathbf{\varepsilon_0}$ is given by
\begin{equation}
    \mathbf{\varepsilon_0}=\phi\beta \quad \text{plane stress},
\end{equation}
\begin{equation}
    \mathbf{\varepsilon_0}=(1+v) \phi \beta \quad \text{plane strain},
\end{equation}
with
\begin{equation}
    \beta=[\alpha \quad \alpha \quad 0]^{\mathrm{T}},
\end{equation}
where $\alpha$ and $v$ denote the coefficient of thermal expansion and Poisson's ratio, respectively. The thermal stress is incorporated as the initial stress $\mathbf{\sigma_0}$, expressed as:
\begin{equation}
   \mathbf{\sigma_0}=\mathbf{D}\mathbf{\varepsilon_0}.
\end{equation}

When calculating thermal stresses, external loads such as body and surface forces are ignored. If present, the effects of these loads, including stress, can be superimposed using the principle of superposition. In the absence of such forces, the total potential energy of an elastic body is solely due to strain energy, given by
\begin{equation}
    \label{eq1}
    \begin{split}
    &\varPi =\frac{1}{2}\int_{v}(\varepsilon^e)^\mathrm{T}\mathbf{\sigma}dv \\
    &=\frac{1}{2}\int_{v}(\varepsilon-\varepsilon_0)^\mathrm{T}\mathbf{D}(\mathbf{\varepsilon-\varepsilon_0})dv \\
    &=\frac{1}{2}\int_{v}\varepsilon^\mathrm{T}\mathbf{D}\mathbf{\varepsilon}dv-\int_{v}\varepsilon^\mathrm{T}\mathbf{D}\mathbf{\varepsilon_0}dv+\frac{1}{2}\int_{v}{\varepsilon_0}^\mathrm{T}\mathbf{D}\mathbf{\varepsilon_0}dv.\\
    \end{split}
\end{equation}

After discretization, and considering that the third term does not affect the variation, it can be omitted. Thus, the above equation simplifies to:
\begin{equation}
    \label{eq2}
    \begin{aligned}
   & \varPi=\frac{1}{2}\sum \limits_{e}(\mathbf{a}^e)^\mathrm{T} \int_{v^e}\mathbf{B}_\mathrm{el}^\mathrm{T}\mathbf{D}\mathbf{B}_\mathrm{el}dv\mathbf{a}^e-\sum\limits_{e}(\mathbf{a}^\mathrm{T})\int_{v^e}\mathbf{B}_\mathrm{el}^\mathrm{T}\mathbf{D}\mathbf{\varepsilon_0}dv \\
   & =\frac{1}{2}\mathbf{a}^\mathrm{T}\mathbf{K}_\mathrm{el}\mathbf{a}-\mathbf{a}^\mathrm{T}\mathbf{f}_\mathrm{el},
    \end{aligned}
\end{equation}
where the matrix $\mathbf{K}_{\mathrm{el}}$ is the global stiffness matrix, given by:
\begin{equation}
    \mathbf{K}_\mathrm{el} = \int_{v^e} \mathbf{B}_{\mathrm{el}}^\mathrm{T} \mathbf{D} \mathbf{B}_\mathrm{el} \, dv.
\end{equation}

Similarly, for the external forces due to initial strains, we can write:
\begin{equation}
    \mathbf{f}_\mathrm{el} = \int_{v^e} \mathbf{B}_\mathrm{el}^\mathrm{T} \mathbf{D} \mathbf{\varepsilon_0} \, dv.
\end{equation}

\subsection{Polygonal element construction}
The Wachspress shape function possesses several important properties, including partition of unity, Kronecker delta function, non-negativity, and linear precision~\cite{nguyen2017polytree}. Its main advantage is the $C^\infty$ smoothness within the polygonal domain and $C^0$ continuity along boundaries~\cite{wu2023polygonal}, which improves solution accuracy and guarantees the convergence of temperature and displacement fields.

\begin{figure}[H]
  \centering
  \includegraphics[width=0.9\textwidth]{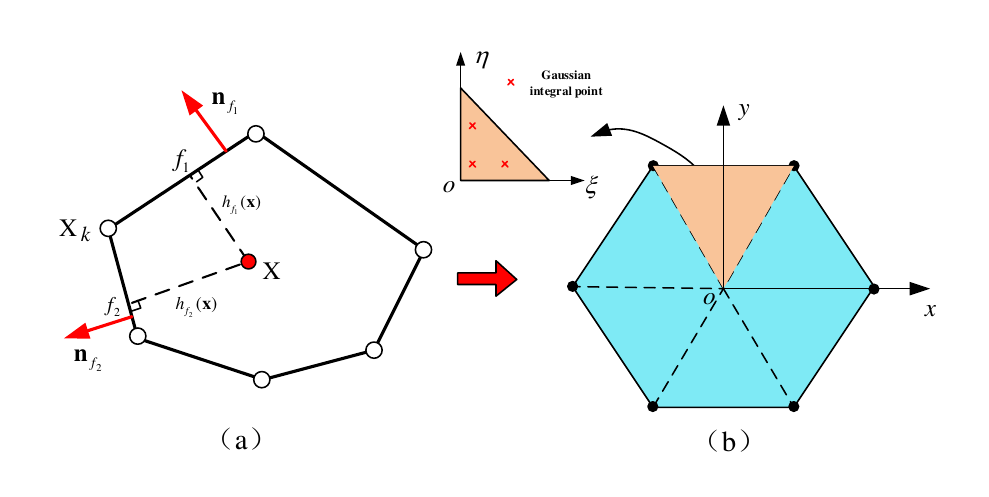}
  \caption{Construction of the polygonal FEM element; (a) Wachspress basis function; (b) numerical integration techniques for polygon.}
  \label{fig:polygonal element construction}
\end{figure}

In this work, we utilize the Wachspress shape function to construct conforming approximations on polygonal meshes. Wachspress introduced rational basis functions for polygonal elements, drawing on projective geometry principles. These functions ensure accurate nodal interpolation and maintain linearity along the boundaries by utilizing algebraic equations for the edges. The polygonal element in Barycentric coordinates is shown in Fig. \ref{fig:polygonal element construction} (a). The shape function $w_k(\mathbf{x})$ is defined as
\begin{equation}
w_k(\mathbf{x})=\frac{\operatorname{det}\left(\mathbf{n}_{f_1}, \mathbf{n}_{f_2}\right)}{h_{f_1}(\mathbf{x}) h_{f_2}(\mathbf{x})},
\end{equation}
where $\mathbf{n}_{f_1}$ and $\mathbf{n}_{f_2}$ denote the normal vectors of the boundary edges. $h_{f_1}(\mathbf{x})$ and $h_{f_2}(\mathbf{x})$ denote the perpendicular distances from $\mathbf{x}_k$ to the facets $f_1$ and $f_2$, which are calculated as follows:
\begin{equation}
h_{f_i}(\mathbf{x}) = (\mathbf{x}_k - \mathbf{x}) \cdot \mathbf{n}_{f_i}.
\end{equation}

The edges $f_i$ should be arranged in a counterclockwise order around the node $\mathbf{x}_k$, as viewed from the exterior. The shape functions $N_k(\mathbf{x})$ are then given by \cite{warren2003uniqueness,warren2007barycentric}
\begin{equation}N_k(\mathbf{x})=\frac{w_\mathbf{k}(\mathbf{x})}{\sum_{\mathbf{k}\in V}w_\mathbf{k}(\mathbf{x})}.\end{equation}

\subsection{Numerical integration technique for polygonal elements}
To compute the stiffness matrix and mass matrix of polygonal elements, this work adopts a triangulation-based numerical integration strategy. Specifically, the polygonal element is first subdivided into several non-overlapping triangles. The integration is then performed using Gaussian quadrature over each triangular subdomain, and the results are assembled to form the stiffness matrix of the entire polygonal element. 

Consider a polygonal element $\Omega$ with $m$ edges. The domain can be subdivided into $m$ triangular subdomains $\Omega_t$ using either an internal node or by connecting all edges to the centroid $\mathbf{x}_c$:
\begin{equation}
\Omega = \bigcup_{t=1}^{m} \Omega_t,
\end{equation}
where $\Omega_t$ denotes the $t$-th triangular subdomain.

For each triangular subdomain $\Omega_t$, Gaussian quadrature is applied:
\begin{equation}
\int_{\Omega_t} f(\mathbf{x}) d\mathbf{x} = \sum_{i=1}^{n_g} w_i f(\mathbf{x}_i),
\end{equation}
where $\mathbf{x}_i$ are the Gaussian integration points within the triangle $\Omega_t$, $w_i$ are the corresponding integration weights, and $n_g$ is the number of quadrature points.

Given the shape functions $N_k(\mathbf{x})$ of the polygonal element, using Gaussian quadrature over triangular subdomains, the final expression of stiffness matrix and mass matrix becomes:
\begin{equation}
\mathbf{K}_{\mathrm{el}}^{e} = \sum_{t=1}^{m} \sum_{i=1}^{n_g} w_i \mathbf{B}_{\mathrm{el}}^{\mathrm{T}}(\mathbf{x}_i) \mathbf{D} \mathbf{B}_{\mathrm{el}}(\mathbf{x}_i) |\mathbf{J}_t|,
\end{equation}
\begin{equation}
    \mathbf{K}_{\mathrm{th}}^{e} = \sum_{t=1}^{m} \sum_{i=1}^{n_g} w_i k\mathbf{B}_{\mathrm{th}}^{\mathrm{T}}(\mathbf{x}_i) \mathbf{B}_{\mathrm{th}}(\mathbf{x}_i) |\mathbf{J}_t|,
\end{equation}
\begin{equation}
    \mathbf{M}_{\mathrm{th}}^{e} = \sum_{t=1}^{m} \sum_{i=1}^{n_g} w_i \rho c\mathbf{N}^{\mathrm{T}}(\mathbf{x}_i) \mathbf{N}(\mathbf{x}_i) |\mathbf{J}_t|.
\end{equation}

\section{Implementation}
\label{sec:implementation}
\subsection{Steady-state solutions}
In this study, we only consider the effect of the temperature field on the displacement field, neglecting the influence of the displacement field on the temperature field. Hence, the thermal stress analysis focuses on the coupling between thermal and stress fields. The coupling equation can be expressed as
\begin{equation}
\begin{bmatrix}
\mathbf{K}_{\mathrm{th}} & \mathrm{0} \\
\mathrm{C}_{el} & \mathbf{K}_{\mathrm{el}}
\end{bmatrix}
\left\{
\begin{array}{c}
\boldsymbol{\phi} \\
\boldsymbol{u}
\end{array}
\right\}
=
\left\{
\begin{array}{c}
\boldsymbol{f}_{\mathrm{th}} \\
\boldsymbol{f}_{\mathrm{el}}
\end{array}
\right\},
\label{eq:th-el-matrix}
\end{equation}
where $\mathbf{K}^{\mathrm{th}}$ and $\mathbf{K}^{\mathrm{el}}$ denote the thermal and elastic stiffness matrices, respectively. The solution vector $\boldsymbol{\phi}$ and $\boldsymbol{u}$ are the temperature field and the displacement field, respectively. The vectors $\boldsymbol{f}^{\mathrm{th}}$ and $\boldsymbol{f}^{\mathrm{el}}$ denote the right-hand of the flux and force, respectively.

Eq. (\ref{eq:th-el-matrix}) can be solved as two steps. Firstly, we solves the temperature field as the following:
\begin{equation}
    \mathbf{K}_{\mathrm{th}}\boldsymbol{\phi}=\boldsymbol{f}_{\mathrm{th}}.
\end{equation}

After obtaining the temperature field, the displacement field is then solved using the following equation:
\begin{equation}
   \mathbf{K}_{\mathrm{el}}\boldsymbol{u}=\boldsymbol{f}_{\mathrm{el}}-\mathrm{C}_\mathrm{el}\boldsymbol{\phi}.
\end{equation}

\subsection{Transient solutions}
The coupling equation of transient thermal stress can be expressed as:
\begin{equation}
\begin{bmatrix}
\mathbf{M}_{\mathrm{th}} & \mathrm{0} \\
\mathrm{0}  & \mathrm{0}
\end{bmatrix}
\left\{
\begin{array}{c}
\dot{\boldsymbol{\phi}} \\
\boldsymbol{u}
\end{array}
\right\}
+
\begin{bmatrix}
\mathbf{K}_{\mathrm{th}} & \mathrm{0} \\
\mathbf{C}_{el} & \mathbf{K}_{\mathrm{el}}
\end{bmatrix}
\left\{
\begin{array}{c}
\boldsymbol{\phi} \\
\boldsymbol{u}
\end{array}
\right\}
=
\left\{
\begin{array}{c}
\boldsymbol{f}_{\mathrm{th}} \\
\boldsymbol{f}_{\mathrm{el}}
\end{array}
\right\},
\label{eq:th-el-transient}
\end{equation}
where $\mathbf{M}_{\mathrm{th}}$ is the mass matrix of heat conduction, $\dot{\boldsymbol{\phi}}$ is temperature change rate.  At time $[t,t+\Delta t]$, the temperature change rate $\dot{\boldsymbol{\phi}}$ can be expressed as 
\begin{equation}
    \boldsymbol{\dot{\phi}}(t)= \frac{\Delta \boldsymbol{\phi}}{\Delta t} = \frac{\boldsymbol{\phi}(t + \Delta t) - \boldsymbol{\phi}(t)}{\Delta t}.
\end{equation}

Similarly, Eq. (\ref{eq:th-el-transient}) can also be solved in two steps: Firstly, we solves the temperature field $\boldsymbol{\phi}$ of each moment as the following:
\begin{equation}
\mathbf{M}_{\mathrm{th}}\boldsymbol{\dot{\phi}}+\mathbf{K}_{\mathrm{th}}\boldsymbol{\phi}=\boldsymbol{f}_{\mathrm{th}}.
\end{equation}

Once the temperature field $\boldsymbol{\phi}$ is determined, the displacement field is subsequently solved using the following equation:
\begin{equation}
   \mathbf{K}_{\mathrm{el}}\boldsymbol{u}=\boldsymbol{f}_{\mathrm{el}}-\mathbf{C}_{el}\boldsymbol{\phi}.
\end{equation}
 
\subsection{Implementation using  ABAQUS UEL}
In this work, we implemented thermal-stress analysis using ABAQUS with User-Defined Elements (UEL). Algorithm \ref{alg:1} presents a flowchart outlining the procedure for both steady-state and transient heat conduction analysis. The primary function of the UEL in ABAQUS is to update the element's contribution to the residual force vector (RHS) and the stiffness matrix (AMATRX) through the user subroutine interface provided by the software. For steady state thermal stress analysis, the AMATRX and RHS are defined as follows:
\begin{equation}
\mathrm{AMATRX}=\begin{bmatrix}
\mathbf{K}_{\mathrm{th}} & \mathrm{0} \\
\mathbf{C}_{el} & \mathbf{K}_{\mathrm{el}}
\end{bmatrix},\label{eq:amatrx-steady-state}
\end{equation}

\begin{equation}
\mathrm{RHS}=
-\begin{bmatrix}
 \mathbf{K}_{\mathrm{th}} & \mathrm{0} \\
\mathrm{C}_{el} & \mathbf{K}_{\mathrm{el}}
\end{bmatrix}
\left\{
\begin{array}{c}
\boldsymbol{\phi} \\
\boldsymbol{u}
\end{array}
\right\}+
\left\{
\begin{array}{c}
\boldsymbol{f}_{\mathrm{th}} \\
\boldsymbol{f}_{\mathrm{el}}
\end{array}
\right\}.\label{eq:rhs-steady-state}
\end{equation}

For transient thermal-stress analysis, the AMATRX and RHS are defined as follows:
\begin{equation}
\mathrm{AMATRX}=\begin{bmatrix}
\mathbf{K}_{\mathrm{th}} & \mathrm{0} \\
\mathbf{C}_{el} & \mathbf{K}_{\mathrm{el}}
\end{bmatrix}^{t+\Delta t}+\frac{1}{\Delta t}{\begin{bmatrix}
\mathbf{M}_{\mathrm{th}} & \mathrm{0} \\
\mathrm{0}  & \mathrm{0}
\end{bmatrix}^{t+\Delta t}}, \label{eq:amatrx-tranisent-state}
\end{equation}

\begin{equation}
    \begin{aligned}
        \mathrm{RHS} = & -\begin{bmatrix}
            \mathbf{K}_{\mathrm{th}} & \mathrm{0} \\
            \mathbf{C}_{\mathrm{el}} & \mathbf{K}_{\mathrm{el}}
        \end{bmatrix}^{t+\Delta t}
        \left\{\begin{array}{c}
            \boldsymbol{\phi} \\
            \boldsymbol{u}
        \end{array}\right\}^{t+\Delta t} \\
        & -\frac{1}{\Delta t} \begin{bmatrix}
            \mathbf{M}_{\mathrm{th}} & \mathrm{0} \\
            \mathrm{0} & \mathrm{0}
        \end{bmatrix}^{t+\Delta t}
        \left(
              \left\{\begin{array}{c}
                \boldsymbol{\phi} \\
                \boldsymbol{u}
            \end{array}\right\}^{t+\Delta t}
            -
              \left\{\begin{array}{c}
                \boldsymbol{\phi} \\
                \boldsymbol{u}
            \end{array}\right\}^{t}
        \right).
    \end{aligned}\label{eq:rhs-tranisent-state}
\end{equation}

\begin{algorithm}[H]
	\caption{Solving the thermal stress problems using the PFEM}\label{alg:1}
    \textbf{Input:} Nodes and elements information, material properties, and nodal temperatures $\boldsymbol{\phi}_t$ and nodal displacements $\boldsymbol{u}_t$ \\
    \textbf{Output:} Nodal temperatures $\boldsymbol{\phi}_{t+1}$ and nodal displacements $\boldsymbol{u}_{t+1}$
	\begin{algorithmic}[1]
        \While{ABAQUS not converged}
        \State{Solve nodal temperatures $\boldsymbol{\phi}^k_{t+1}$ and nodal displacements $\boldsymbol{u}^k_{t+1}$}
		\For {1 to AllEle}  \quad \quad  \textit{! Traverse all elements}
        \For{1 to Element subdomains}
        \State{Construct triangular subdomains based on polygons}
        \State{Solve the subdomains of the stiffness $\mathbf{K}_{\mathrm{el}}^{sub}$ and $\mathbf{K}_{\mathrm{th}}^{\mathrm{sub}}$ and mass matrix $\mathbf{M}_{\mathrm{th}}^{\mathrm{sub}}$}
        \EndFor
         \State{Assemble the element stiffness $\mathbf{K}_{\mathrm{el}}$ and $\mathbf{K}_{\mathrm{th}}$ and mass matrix $\mathbf{M}_{\mathrm{th}}$}
         \If{lflags(1)=71}
         \State{Update stiffness matrix $\mathbf{AMATRX}$ and residual force vector $\mathbf{RHS}$ using Eqs. (\ref{eq:amatrx-steady-state}) and (\ref{eq:rhs-steady-state})}
         \EndIf
         \If{lflags(1)=72 or lflags(1)=73}
         \State{Update stiffness matrix $\mathbf{AMATRX}$ and residual force vector $\mathbf{RHS}$ using Eqs. (\ref{eq:amatrx-tranisent-state}) and (\ref{eq:rhs-tranisent-state})}
         \EndIf
		\EndFor
        \State{\textbf{update} $k=k+1$}
        \State{Solve nodal temperatures $\boldsymbol{\phi}_{t+1}$=$\boldsymbol{\phi}^k_{t+1}$ and nodal displacements $\boldsymbol{u}_{t+1}$=$\boldsymbol{u}^k_{t+1}$}
        \EndWhile
	\end{algorithmic} 
\end{algorithm} 

\subsection{Acceleration technology based on the quadtree parent element}
In numerical analysis, computational models are generally discretized into numerous finite elements, each of which requires individual calculations. Enhancing the computational efficiency at the element level can thus lead to a substantial improvement in the overall performance of the simulation.

The relationship between a parent element and its sub-element can be expressed through the following scaling mappings:
\begin{equation}
    \mathbf{K}_{\mathrm{th}}^{\mathrm{sub}} = k\,\mathbf{K}_{\mathrm{th}}^{\mathrm{par}}, \label{eq:map1}
\end{equation}
\begin{equation}
    \mathbf{K}_{\mathrm{el}}^{\mathrm{sub}} = E\,\mathbf{K}_{\mathrm{el}}^{\mathrm{par}},
\end{equation}
\begin{equation}
    \mathbf{M}_{\mathrm{th}}^{\mathrm{sub}} = \rho\,c\,L\,\mathbf{M}_{\mathrm{th}}^{\mathrm{par}}, \label{eq:map2}
\end{equation}
where $\mathbf{K}^{\mathrm{sub}}$ and $\mathbf{M}^{\mathrm{sub}}$ denote the stiffness and mass matrices of the sub-element, respectively, and $\mathbf{K}^{\mathrm{par}}$ and $\mathbf{M}^{\mathrm{par}}$ represent those of the parent element. The parameters $E$, $k$, $\rho$, $c$, and $L$ correspond to the Young’s modulus, thermal conductivity, density, specific heat capacity, and characteristic length of the sub-element, respectively.

\begin{figure}[H]
  \centering
  \includegraphics[width=1.0\textwidth]{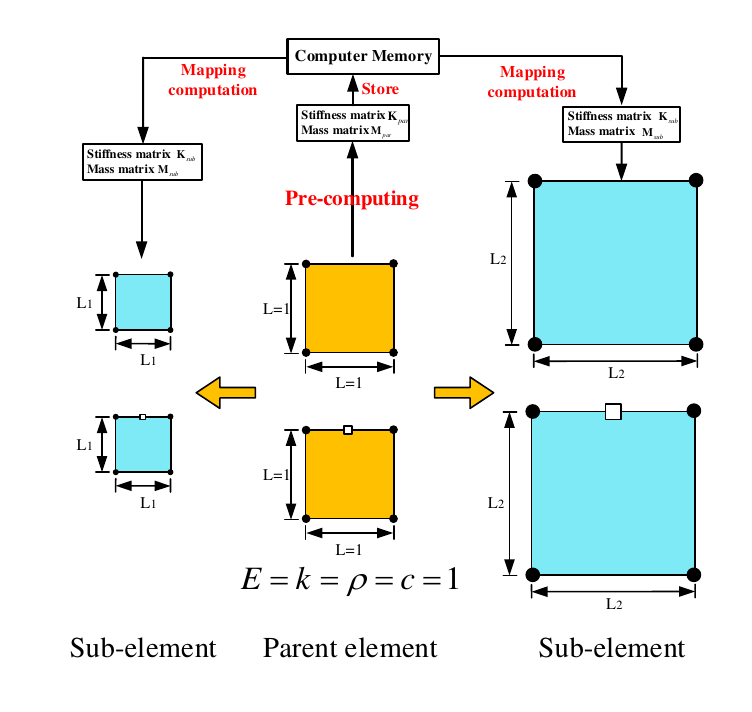}
  \caption{Mapping relationship between the parent element and sub-element.}
  \label{fig:acctech}
\end{figure}

In the proposed framework, we define parent elements whose geometric and material parameters are set to 1, as presented in Fig. \ref{fig:acctech}. The stiffness and mass matrices of the parent element, $\mathbf{K}^{\mathrm{par}}$ and $\mathbf{M}^{\mathrm{par}}$, are precomputed and stored in memory prior to the analysis. For quadtree-based mesh structures, the matrices of sub-elements are efficiently derived using the scaling relationships defined in Eqs.~(\ref{eq:map1})$\sim$(\ref{eq:map2}), which involve only the rescaling of material and geometric properties. This approach eliminates the need for reassembling element matrices from scratch, thereby achieving a significant reduction in computational cost.

\subsection{Defining the UEL elements}
\begin{figure}[H]
  \centering
  \includegraphics[width=1.0\textwidth]{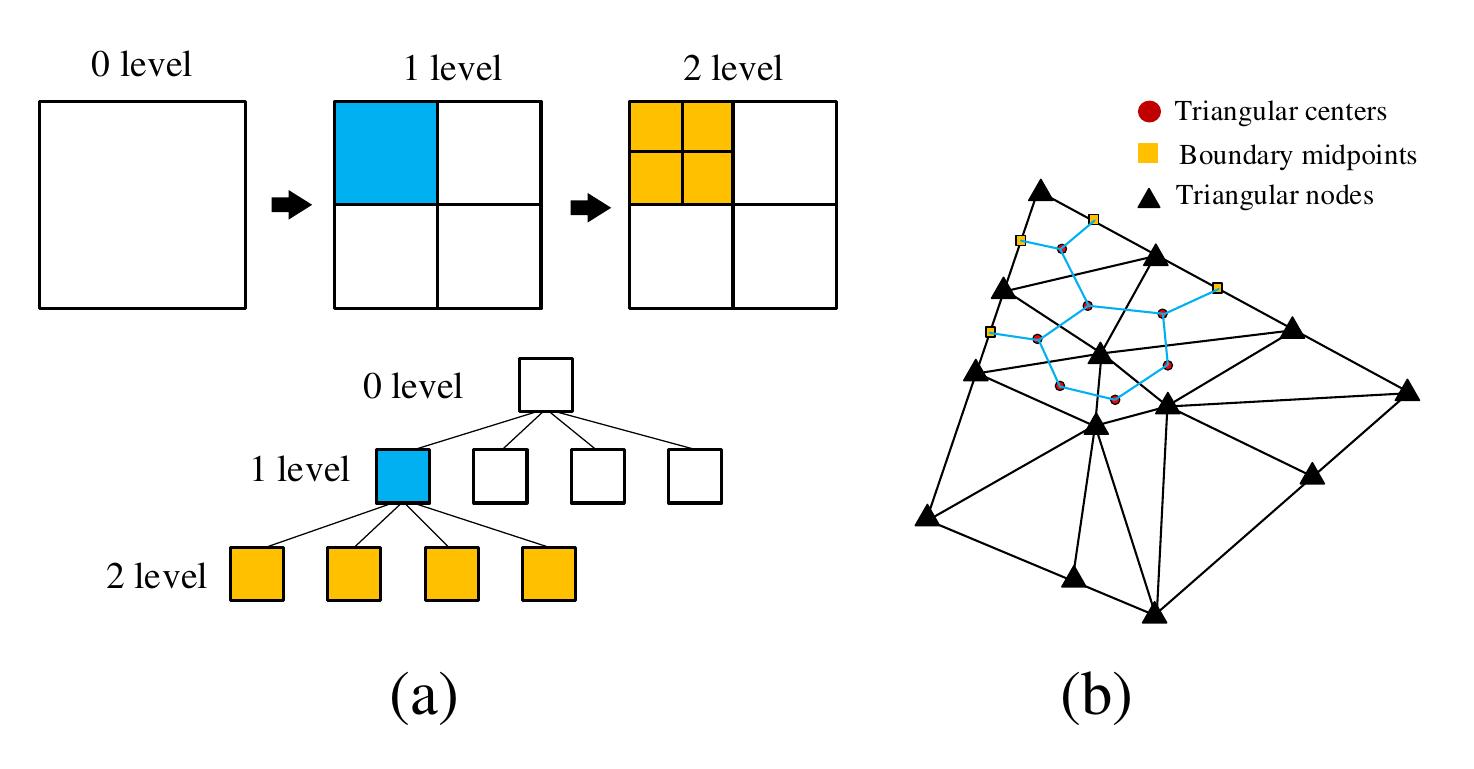}
  \caption{Discretization techniques of the quadtree mesh and polygonal mesh.}
  \label{fig:meshtech}
\end{figure}

In this study, polygonal and quadtree mesh discretization techniques are employed for domain partitioning, with the specific meshing algorithm illustrated in Fig.~\ref{fig:meshtech}. Further details regarding the mesh generation procedures can be found in Refs.~\cite{yePSBFEM2021, yan2024fast}.

\begin{figure}[H]
  \centering
  \includegraphics[width=1.0\textwidth]{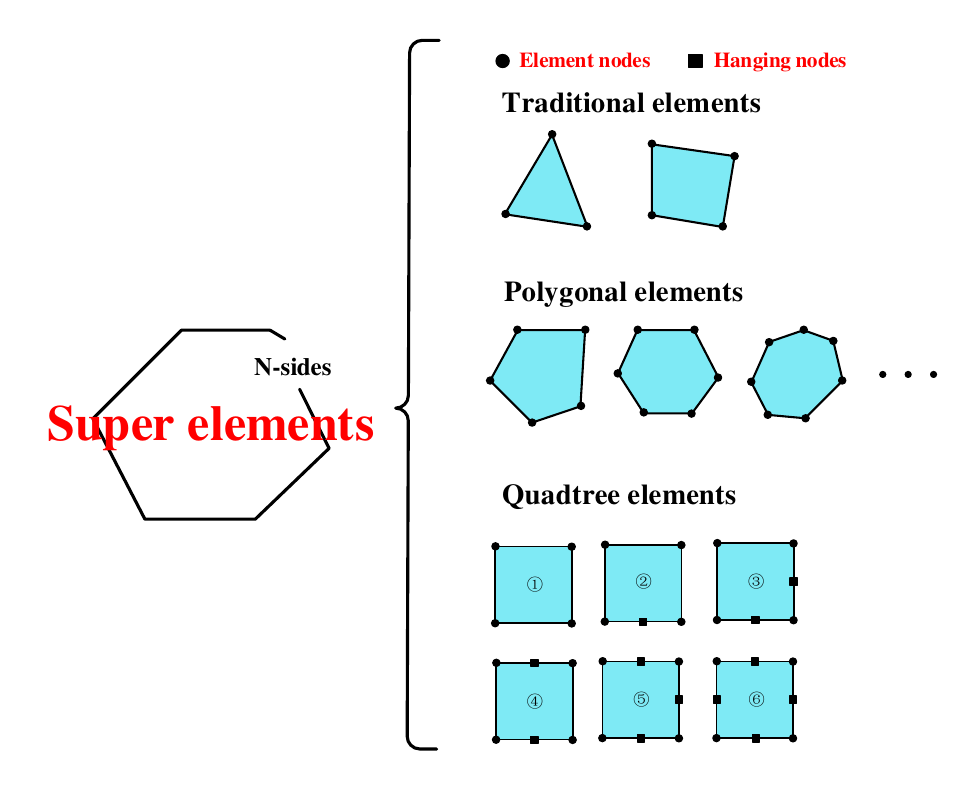}
  \caption{Polygonal finite element method support element types.}
  \label{fig:N_sides_element}
\end{figure}

As depicted in Fig.~\ref{fig:N_sides_element}, a unified polygonal element, termed a “super element”, is developed to accommodate conventional elements, arbitrary polygonal elements, and quadtree-based elements. The numerical model is fully described in the ABAQUS input file, which specifies the nodes, elements, degrees of freedom, and material properties.

An example definition of a hexagonal element (U6) is presented in Listing~\ref{lst:polygonal_input}. Line 1 specifies the element type, number of nodes, element properties, and the degrees of freedom per node. Line 2 defines the active degrees of freedom for displacement and temperature. Line 3 indicates the user-defined element (UEL) type, while Line 4 provides the nodal connectivity. Other element types can be defined analogously by adhering to this structured format.

\begin{lstlisting}[caption={Input file of the polygonal elements in UEL.}, label={lst:polygonal_input}, frame=single, columns=flexible, basicstyle=\ttfamily\small, numbers=left]
*USER ELEMENT,NODES=6,TYPE=U6,PROPERTIES=5,COORDINATES=2,VARIABLES=18
1,2,11
*Element,TYPE=U6
1, 1, 2, 3, 4, 5, 6
*USER ELEMENT,NODES=8,TYPE=U8,PROPERTIES=5,COORDINATES=2,VARIABLES=24
1,2,11
*Element,TYPE=U6
2, 7, 8, 9, 10, 11, 12, 13, 14
\end{lstlisting}

\section{Numerical examples}
\label{sec:examples}
This section presents four numerical benchmarks involving transient thermo-mechanical analysis. The results obtained using the proposed method are compared with those generated by the finite element method (FEM) implemented in the commercial software ABAQUS. To assess the accuracy of the proposed approach, the relative $L_2$ norm of the error, denoted as $\mathbf{e}_{L_2}$, is defined as
\begin{equation}
    \mathbf{e}_{L_2} = \frac{\left\| \mathbf{U}_{\mathrm{num}} - \mathbf{U}_{\mathrm{ref}} \right\|_{L_2}}{\left\| \mathbf{U}_{\mathrm{ref}} \right\|_{L_2}},
\end{equation}
where $\mathbf{U}_{\mathrm{num}}$ denotes the numerical solution, and $\mathbf{U}_{\mathrm{ref}}$ represents the analytical or reference solution. In cases where analytical solutions are not available, the ABAQUS results are used as reference solutions.

\subsection{Steady-state thermal stress problems}
\subsubsection{2D ring plate}
\begin{figure}[H]
  \centering
  \includegraphics[width=0.8\textwidth]{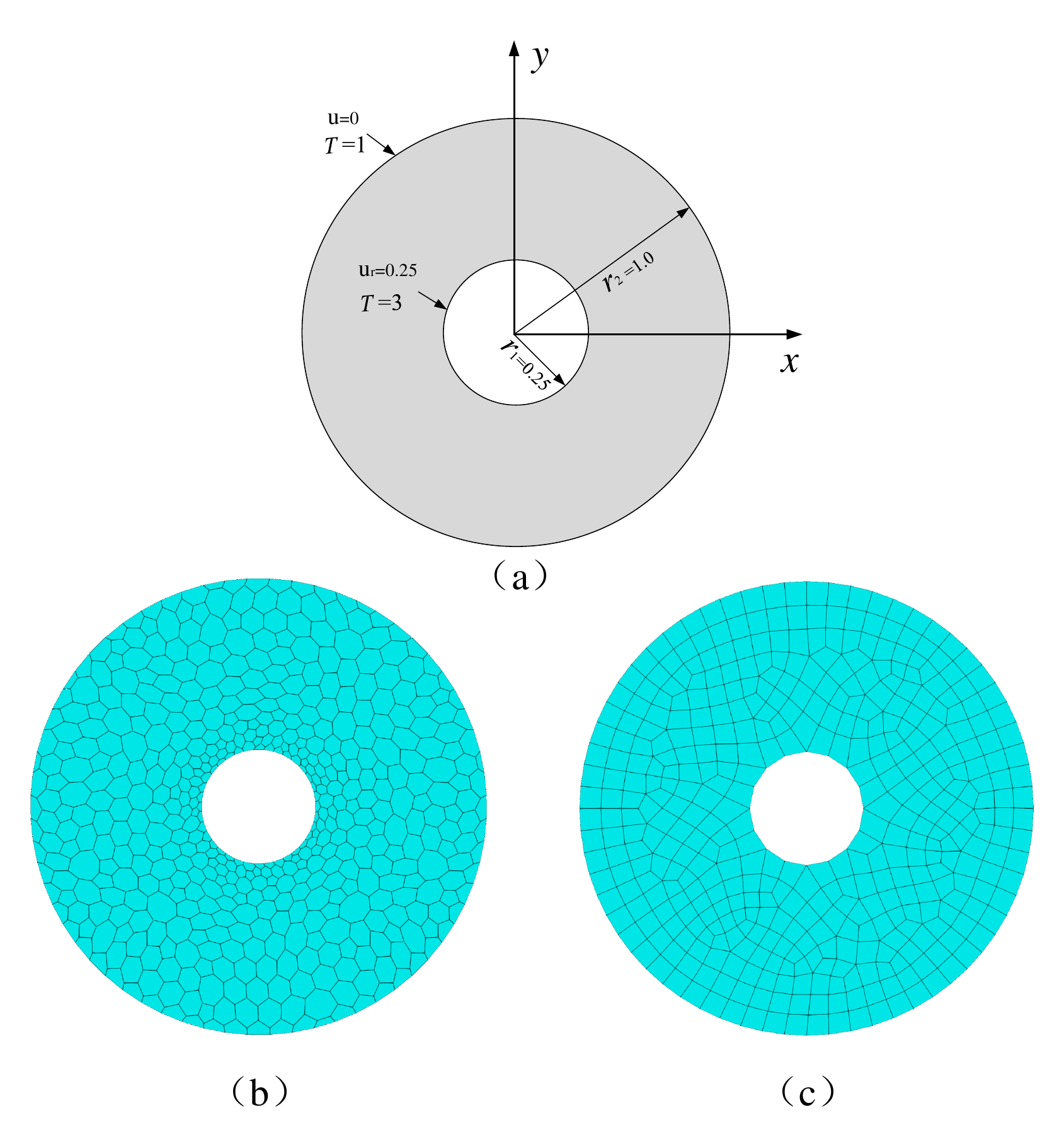}
  \caption{The model of the 2D ring plate; (a) geometry and boundary conditions; (b) polygonal mesh at size 0.1 m; (c) quadrilateral mesh at size 0.1 m.}
  \label{fig:ex02GeoandMesh}
\end{figure}

In the first example, the transient thermoelastic behavior of a two-dimensional ring plate is investigated, as illustrated in Fig.~\ref{fig:ex02GeoandMesh} (a). The ring features an inner radius of 0.25 and an outer radius of 1.0. A radial displacement of 0.25 is prescribed on the inner boundary, while the outer boundary is fixed. Thermal boundary conditions are imposed by applying temperatures of 3 and 1 on the inner and outer boundaries, respectively. The material properties are defined as follows: Young’s modulus $E = 1.0$, Poisson’s ratio $v = 0$, thermal expansion coefficient $\alpha = 1.0$, thermal conductivity $k = 1.0$, and volumetric heat capacity $\rho c = 1.0$.

The analytical solutions for the temperature, radial displacement, and radial stress in cylindrical coordinates are given in Ref.~\cite{li2018solution}:
\begin{equation}
    T_r(r) = 1 - \frac{\ln(r)}{\ln(2)},
\end{equation}
\begin{equation}
    u_r(r) = -\frac{r}{2} \cdot \frac{\ln(r)}{\ln(2)},
\end{equation}
\begin{equation}
    \sigma_r(r) = -1 + \frac{\ln(r) - 1}{2 \ln(2)}.
\end{equation}

\begin{figure}[H]
  \centering
  \includegraphics[width=1.0\textwidth]{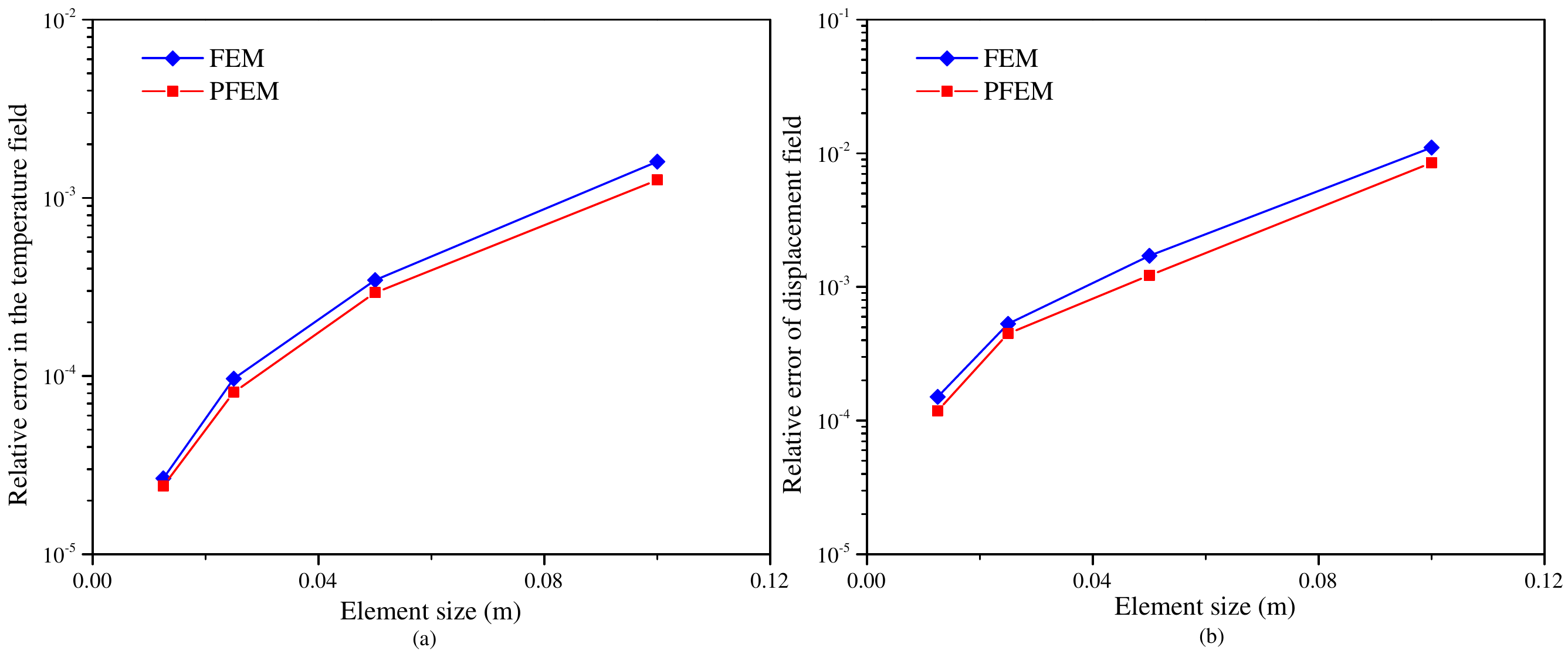}
  \caption{Convergence of the relative errors; (a) temperature field; (b) displacement field.}
  \label{fig:ex02_error}
\end{figure}

\begin{figure}[H]
  \centering
  \includegraphics[width=1.0\textwidth]{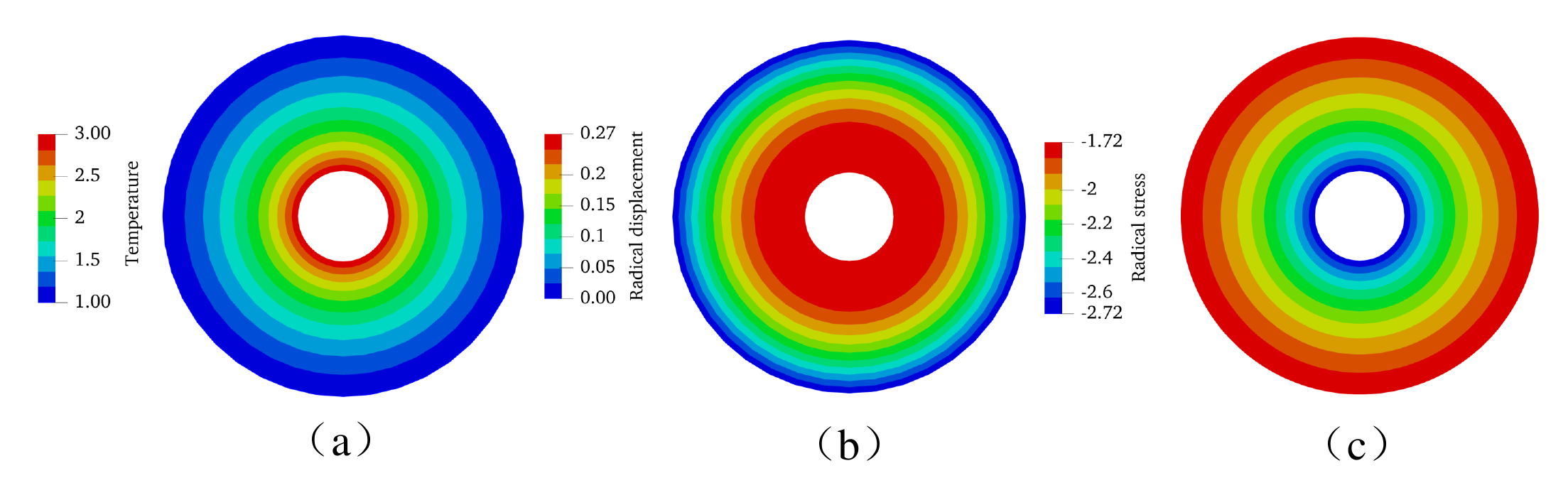}
  \caption{Exact solutions in the radial direction for 2D ring plate; (a) temperature field; (b) radical displacement; (c) radical stress.}
  \label{fig:ex02Geo_analytical}
\end{figure}

\begin{figure}[H]
  \centering
  \includegraphics[width=1.0\textwidth]{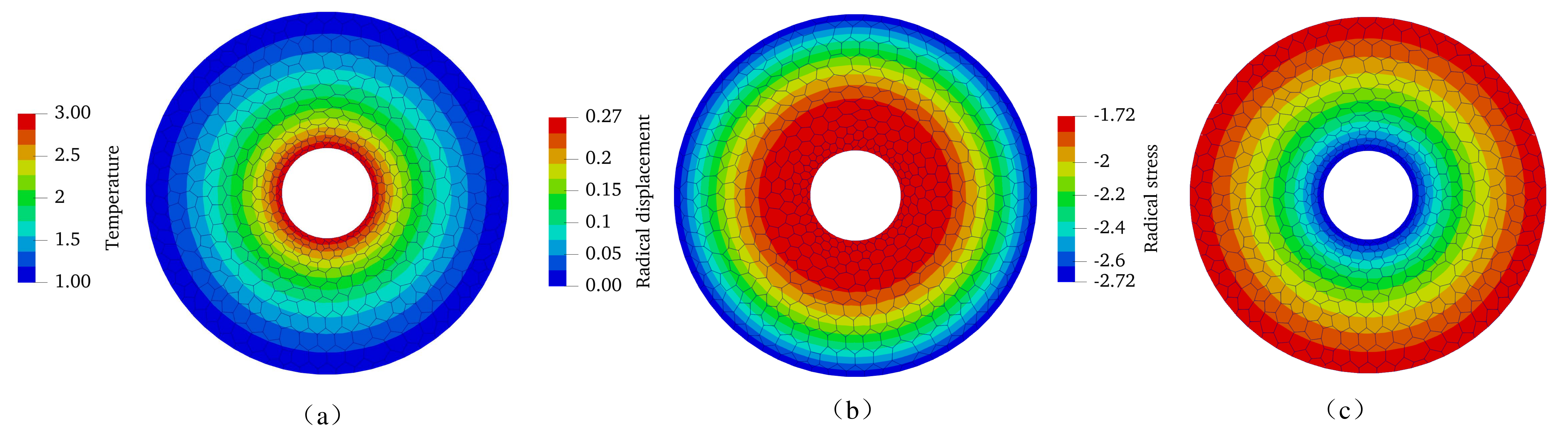}
  \caption{PFEM solutions in the radial direction for 2D ring plate; (a) temperature field; (b) radical displacement; (c) radical stress.}
  \label{fig:ex02_ring01Result}
\end{figure}

To evaluate the accuracy and convergence performance of the proposed method, a sequence of polygonal meshes with element sizes of 0.1 m, 0.05 m, 0.025 m, and 0.0125 m is employed. Figs.~\ref{fig:ex02GeoandMesh} (b) and (c) display the mesh corresponding to an element size of 0.1 m. As shown in Fig.~\ref{fig:ex02_error}, both the proposed PFEM and the conventional FEM demonstrate consistent convergence as the mesh is refined. Notably, PFEM yields superior accuracy compared to FEM for equivalent mesh resolutions. The radial distributions of PFEM results are presented in Fig.~\ref{fig:ex02_ring01Result}, exhibiting excellent agreement with the analytical solutions depicted in Fig.~\ref{fig:ex02Geo_analytical}, thereby validating the effectiveness and accuracy of the proposed formulation.

\subsubsection{A rectangular plate}
\begin{figure}[H]
  \centering
  \includegraphics[width=0.8\textwidth]{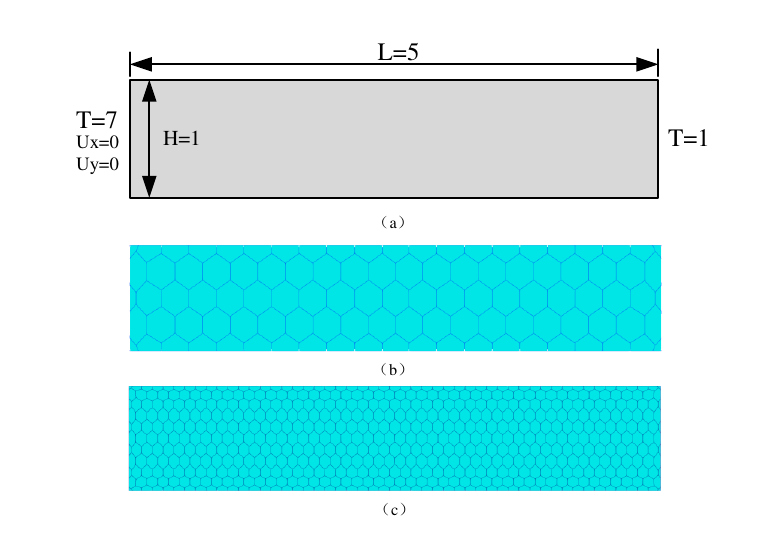}
  \caption{Geometry and mesh model of the rectangular plate; (a) geometry and boundary conditions; (b) 0.25 mesh size; (c) 0.1 mesh size.}
  \label{fig:ex01_geo}
\end{figure}

In this example, a steady-state thermal stress problem of a rectangular plate is analyzed, as illustrated in Fig.~\ref{fig:ex01_geo} (a). The plate has a length of 5 and a height of 1. The left boundary is fully constrained ($\mathrm{U_x = 0}$, $\mathrm{U_y = 0}$), with a temperature of 7 applied, while the right boundary is maintained at a temperature of 1. The material properties are specified as follows: Young’s modulus $E = 1.0$, Poisson’s ratio $v = 0$, thermal expansion coefficient $\alpha = 1.0$, thermal conductivity $k = 1.0$, and volumetric heat capacity $\rho c = 1.0$.

\begin{figure}[H]
  \centering
  \includegraphics[width=0.8\textwidth]{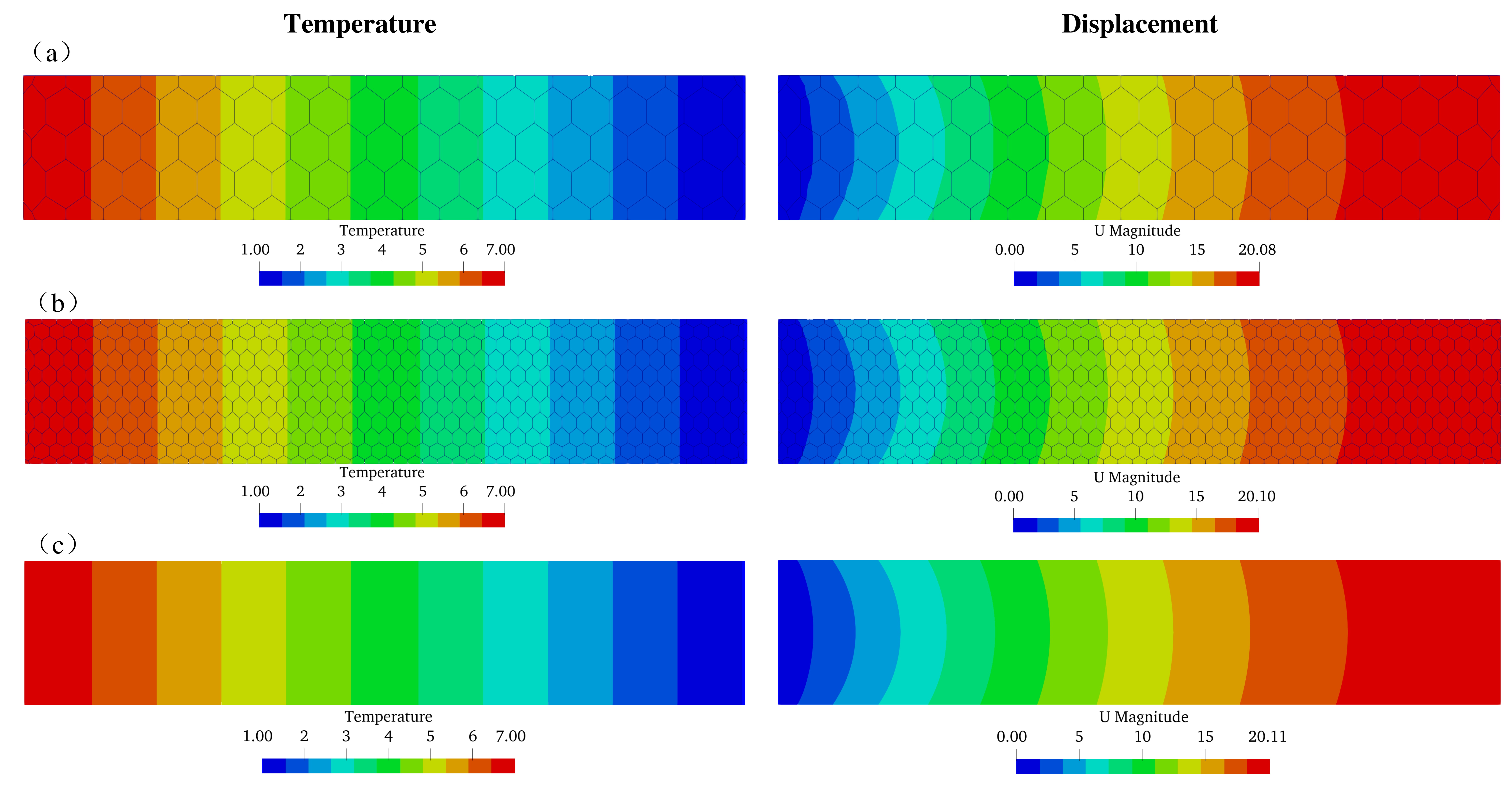}
  \caption{Contours for temperature and displacement fields; (a) 0.25 mesh size; (b) 0.1 mesh size; (c) the reference solution.}
  \label{fig:ex01result}
\end{figure}

The domain is discretized using polygonal elements. A mesh convergence study is performed by refining the element size from 0.5 to 0.25, 0.1, and 0.05. For comparison, the same problem is solved using quadrilateral elements (CPS4T) in ABAQUS. Figs.~\ref{fig:ex01_geo} (b) and (c) show the meshes corresponding to element sizes of 0.25 and 0.1, respectively. Fig.~\ref{fig:ex01result} displays the distributions of temperature and displacement, illustrating excellent agreement between the proposed PFEM and the reference FEM solutions.

\begin{figure}[H]
  \centering
  \includegraphics[width=1.0\textwidth]{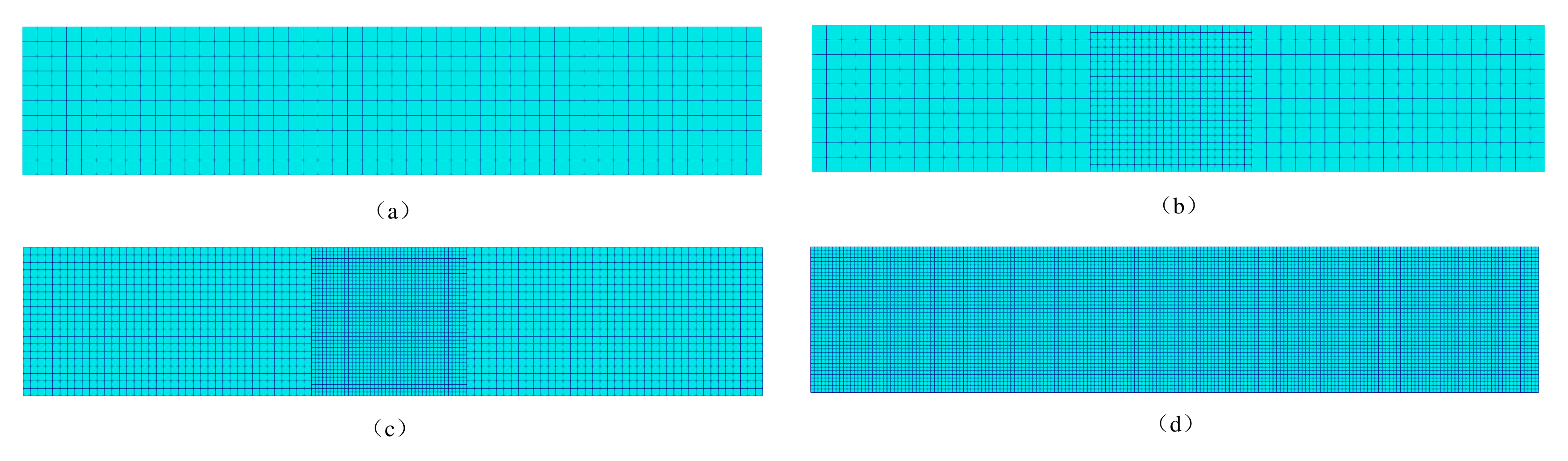}
  \caption{Local mesh refinement model of the beam using the quadtree mesh; (a) coarse mesh; (b) local mesh refinement 1; (c) local mesh refinement 2; (d) fine mesh.}
  \label{ex01_mulscale_mesh}
\end{figure}

\begin{figure}[H]
  \centering
  \includegraphics[width=1.0\textwidth]{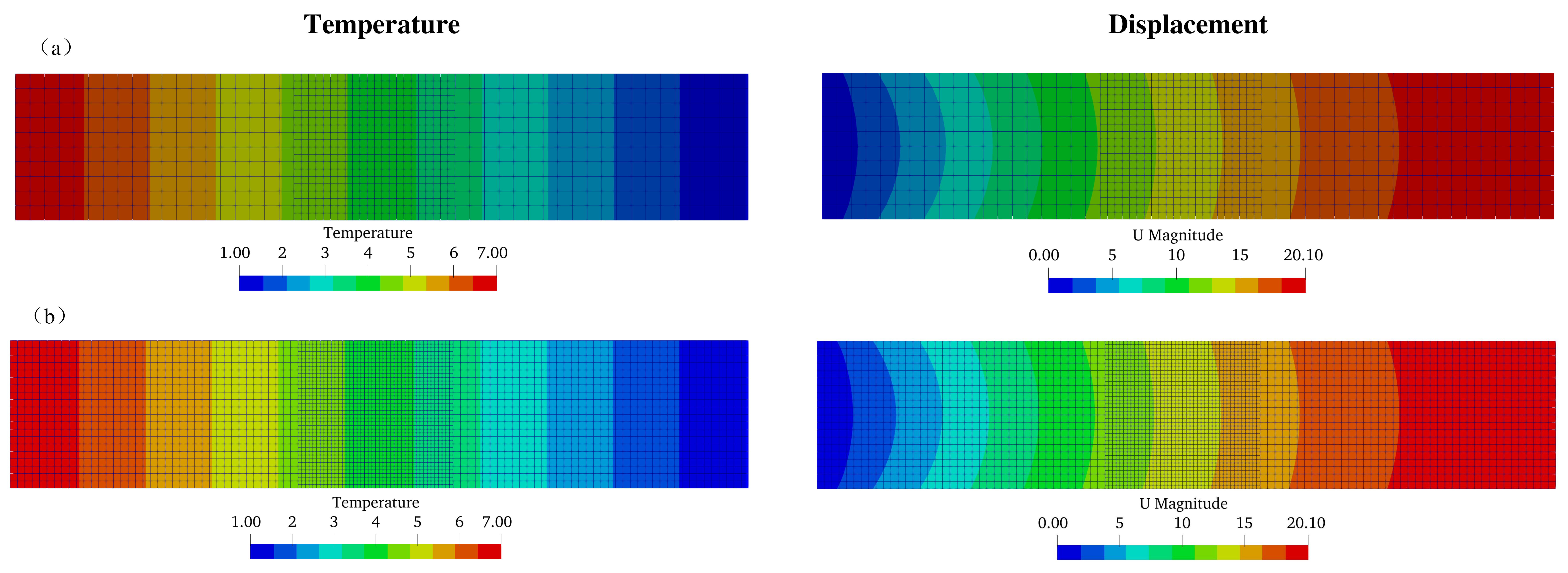}
  \caption{Local mesh refinement model of the beam using the quadtree mesh; (a) coarse mesh; (b) local mesh refinement 1; (c) local mesh refinement 2; (d) fine mesh.}
  \label{ex01_multi-scale_results}
\end{figure}

To further demonstrate the local refinement capability of the proposed method, quadtree meshes are adopted, as shown in Fig.~\ref{ex01_mulscale_mesh}. Figs.~\ref{ex01_mulscale_mesh} (a) and (d) depict coarse and fine global meshes, respectively, while Figs.~\ref{ex01_mulscale_mesh} (b) and (c) present meshes with local refinement focused in the mid-region of the plate. The corresponding mesh statistics and relative errors are summarized in Table~\ref{tab:ex01_t1}. Both the locally refined mesh (Refinement 2) and the globally refined mesh achieve relative errors on the order of $10^{-5}$, confirming the high accuracy of the proposed approach. Notably, the computational time for the locally refined mesh is significantly less than that of the globally refined case, indicating that local refinement can effectively balance computational cost and accuracy. Furthermore, the temperature field shown in Fig.~\ref{ex01_multi-scale_results} confirms that the quadtree-mesh-based solution closely aligns with the reference results.

\begin{figure}[H]
  \centering
  \includegraphics[width=1.0\textwidth]{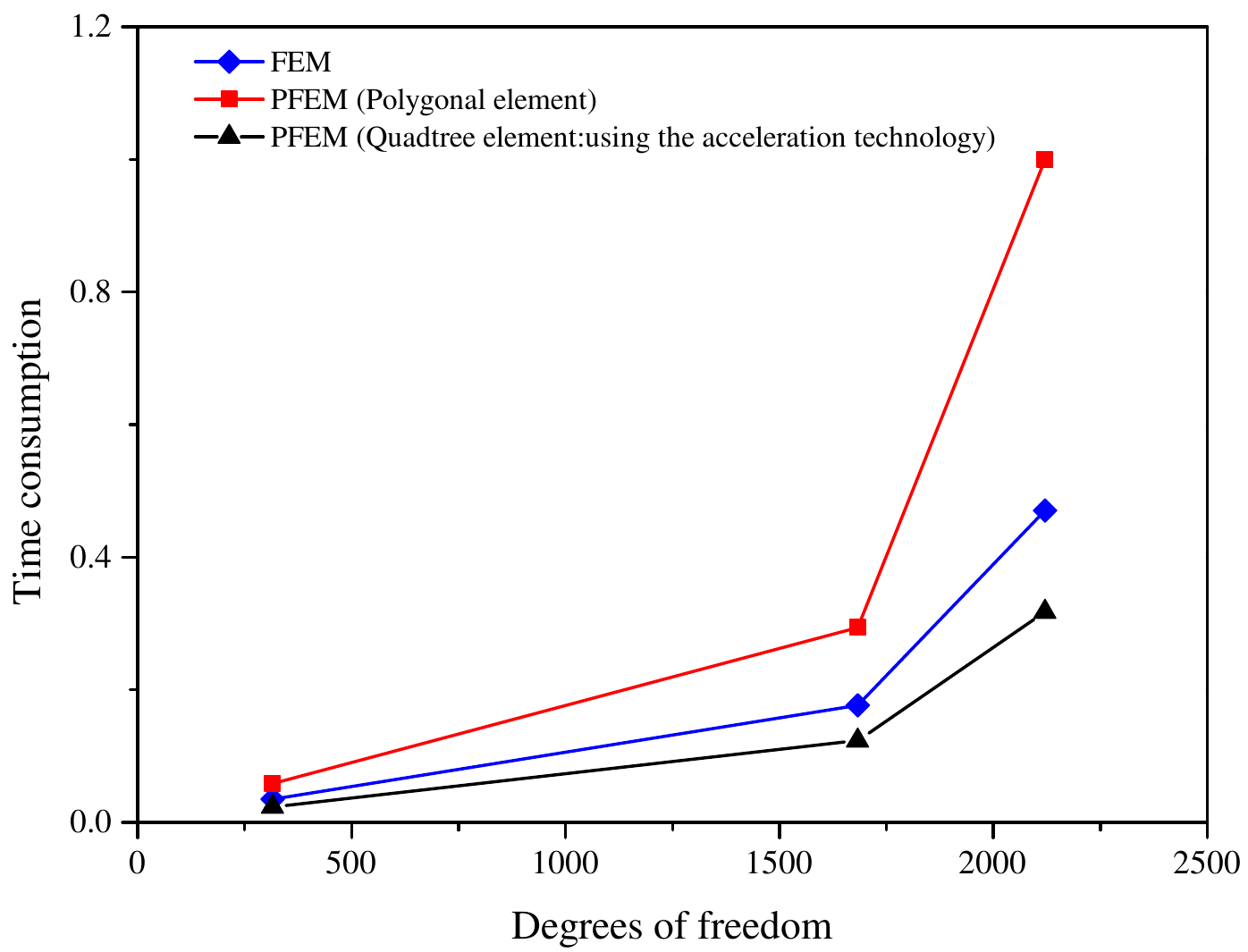}
  \caption{Time consumption comparison for the steady-state thermal stress analysis.}
  \label{fig:ex01time}
\end{figure}

A comparison of computational times for the steady-state thermal stress analysis is presented in Fig.~\ref{fig:ex01time}. Although PFEM generally incurs a slightly higher computational cost than FEM—approximately 1.87 times on average—this increase is attributed to the additional numerical integration required for polygonal elements. However, when employing acceleration techniques, the computational time of PFEM with quadtree meshes is reduced to approximately 44\% of that of the FEM, thereby demonstrating the proposed method’s significant advantage in computational efficiency.

\begin{table}[H]
\centering
\normalsize
\caption{Mesh characteristics and relative errors for the local mesh refinement.}
\resizebox{\textwidth}{!}{ 
\begin{tabular}{cccccc}
\toprule
Mesh type  & Elements & Nodes & Relative error &  CPU time (s) & CPU time (s)*\\ 
\midrule
Coarse mesh            & 500       & 561        & 8.03$\times10^{-4}$  & 0.40  & 0.20\\ 
Fine mesh              & 8000      & 8241      & 1.46$\times10^{-5}$   & 6.70  & 2.60 \\ 
Local mesh refinement 1     & 810      & 912       & 1.09$\times10^{-4}$   & 0.60  & 0.30 \\ 
Local mesh refinement 2     & 3220      & 3422       & 2.19$\times10^{-5}$  & 2.80 &1.10  \\ 
\bottomrule 
\label{tab:ex01_t1}
\end{tabular}}
\vspace{-2em}
\begin{flushleft}
\small Note: * denotes the use of acceleration technology based on the parent element.
\end{flushleft}
\end{table}

\subsection{Transient thermal stress problems}
\subsubsection{L-shaped thin plate}
\begin{figure}[H]
  \centering
  \includegraphics[width=0.8\textwidth]{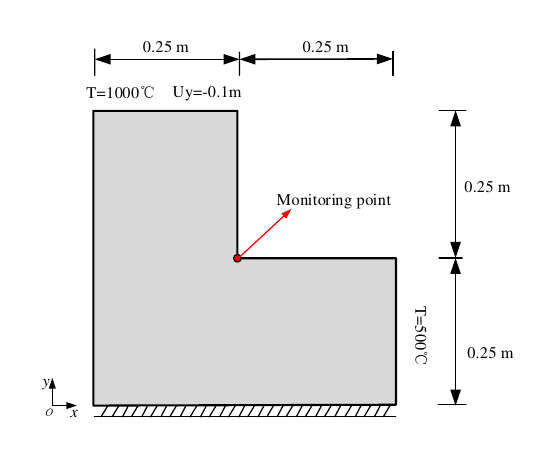}
  \caption{Schematic diagram of a square body with multiple holes.}
  \label{fig:ex03_geo}
\end{figure}

In this example, an L-shaped thin plate is analyzed, as shown in Fig.~\ref{fig:ex03_geo}. The bottom edge of the plate is fully fixed, and a prescribed vertical displacement of \(0.1\,\text{m}\) is applied to the top boundary. Thermal boundary conditions are imposed with the right edge maintained at \(500^\circ\text{C}\) and the top edge at \(1000^\circ\text{C}\). The material properties used in this analysis are: thermal conductivity \(k = 3\,\text{W/m·}^\circ\text{C}\), density \(\rho = 2000\,\text{kg/m}^3\), Young’s modulus \(E = 1 \times 10^4\,\text{Pa}\), Poisson’s ratio \(\nu = 0.3\), thermal expansion coefficient \(\alpha = 0.0011\,^\circ\text{C}^{-1}\), and specific heat \(c = 0.45\,\text{J/g·}^\circ\text{C}\).

\begin{figure}[H]
  \centering
  \includegraphics[width=0.8\textwidth]{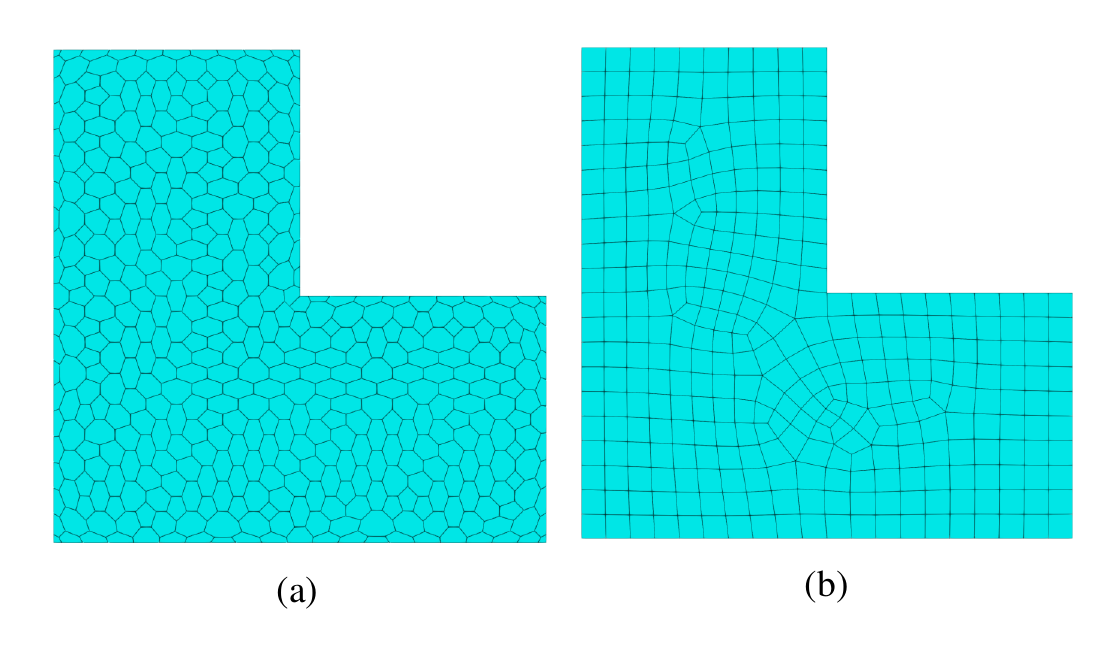}
  \caption{Meshes of the L-shaped thin plate; (a) polygonal mesh; (b) quadrilateral mesh.}
  \label{fig:ex03_mesh}
\end{figure}

To assess the accuracy and convergence of the proposed PFEM for transient thermo-mechanical problems, the domain is discretized using polygonal meshes with element sizes of 0.05~m, 0.025~m, and 0.001~m. A representative mesh is shown in Fig.~\ref{fig:ex03_mesh}. The total simulation time is set to 100~s, with a uniform time step of 1~s.

\begin{figure}[H]
  \centering
  \includegraphics[width=1.0\textwidth]{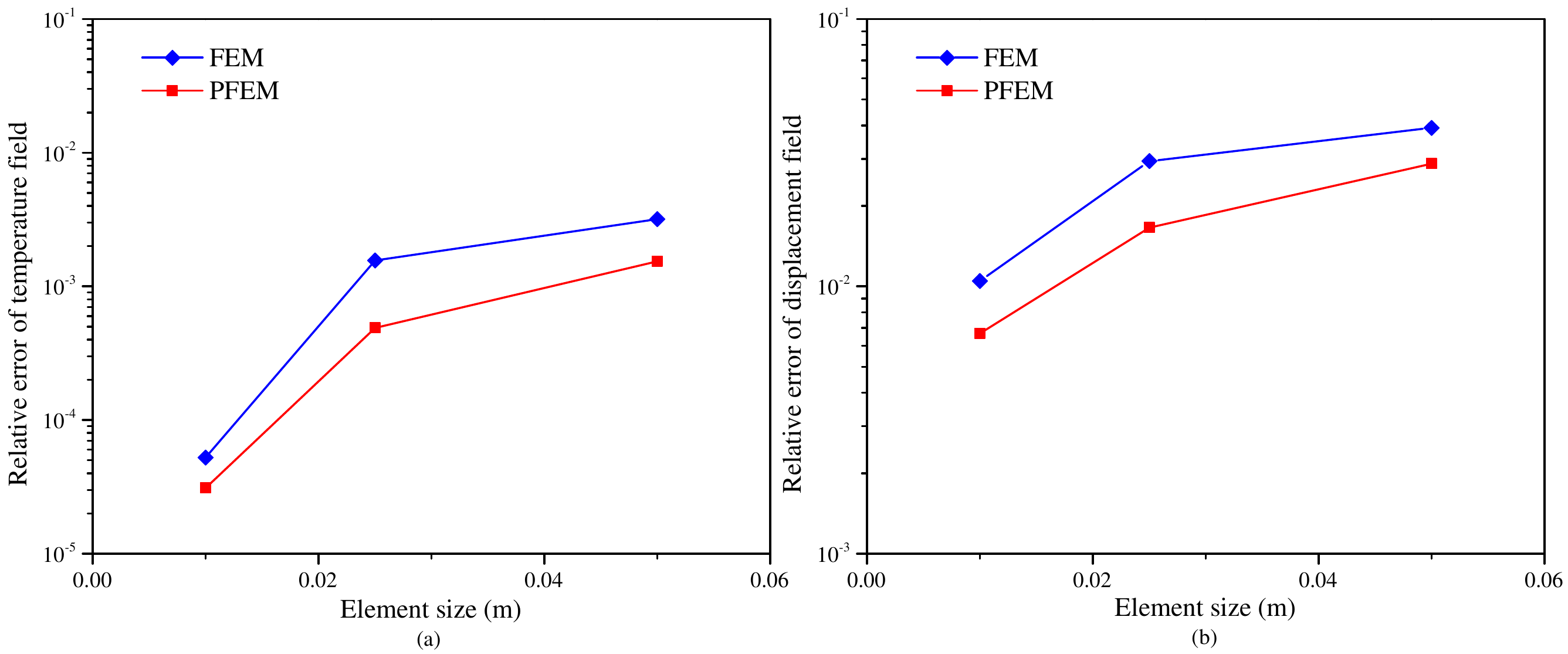}
  \caption{Convergence of the relative errors; (a) temperature field; (b) displacement field.}
  \label{fig:ex03_error}
\end{figure}

\begin{figure}[H]
  \centering
  \includegraphics[width=0.8\textwidth]{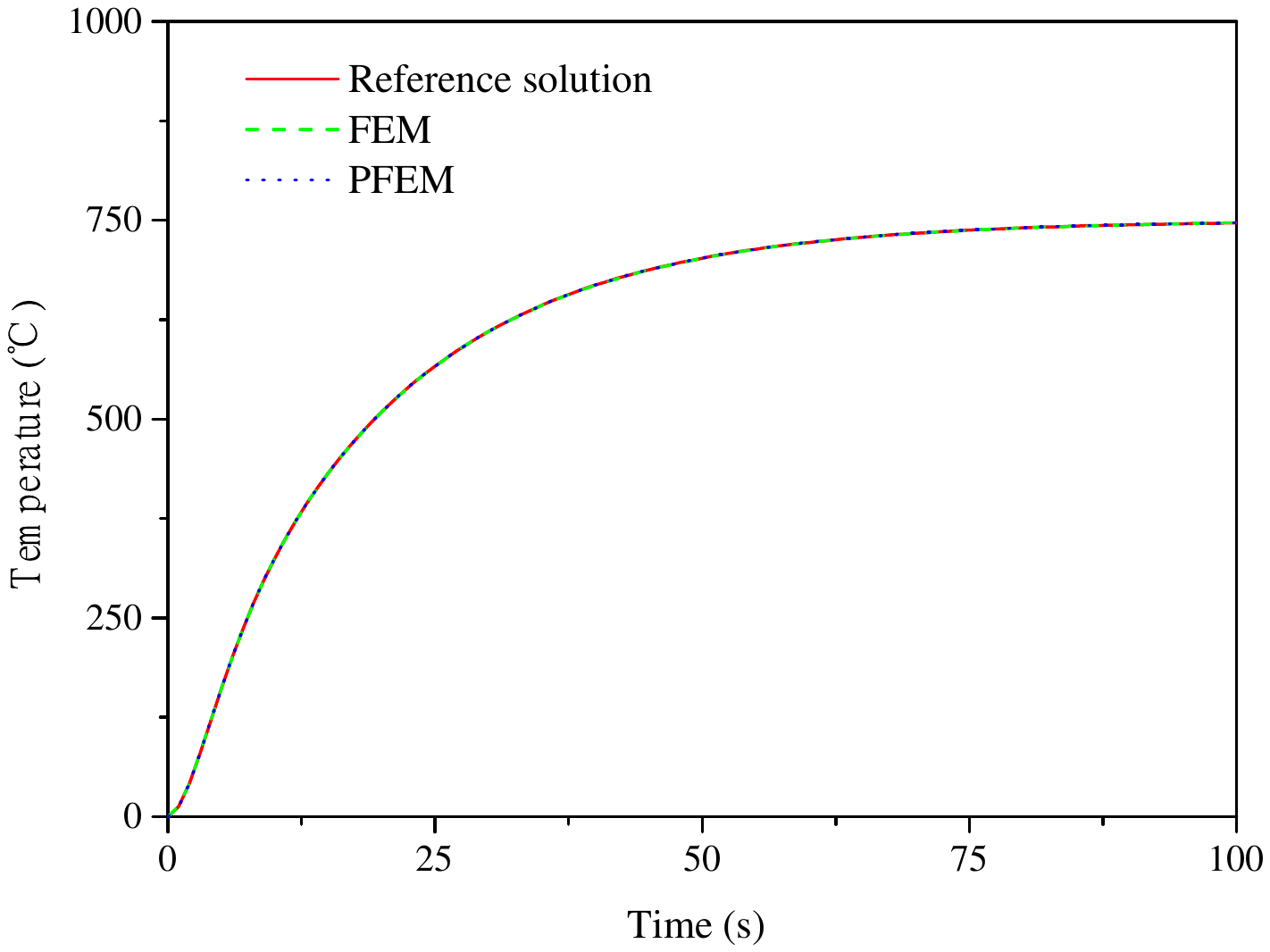}
  \caption{Comparison of the PFEM and FEM in the history of temperature.}
  \label{fig:ex03_temp}
\end{figure}

\begin{figure}[H]
  \centering
  \includegraphics[width=0.8\textwidth]{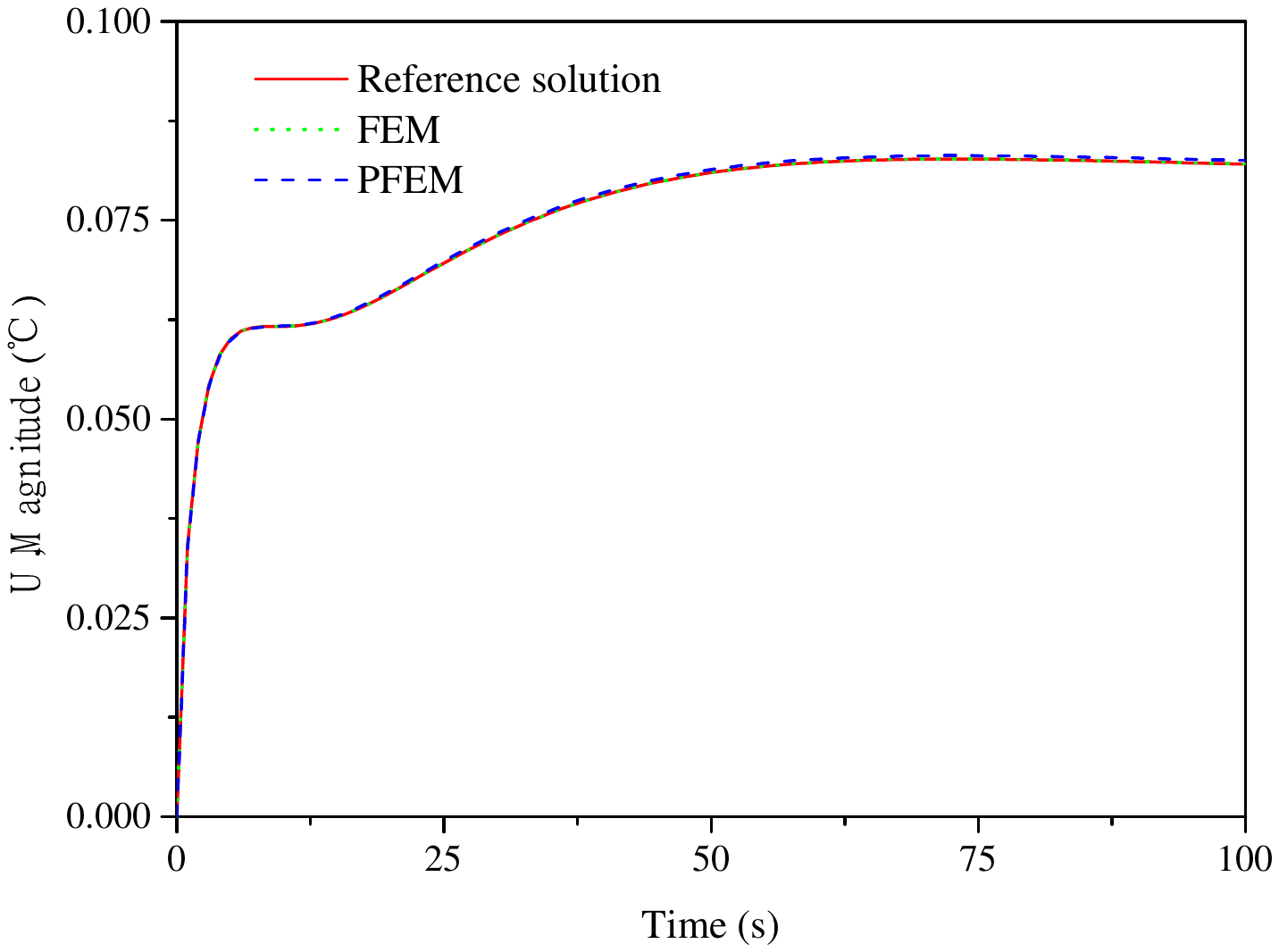}
  \caption{Comparison of the PFEM and FEM in the history of displacement.}
  \label{fig:ex03_disp}
\end{figure}

\begin{figure}[H]
  \centering
  \includegraphics[width=1.0\textwidth]{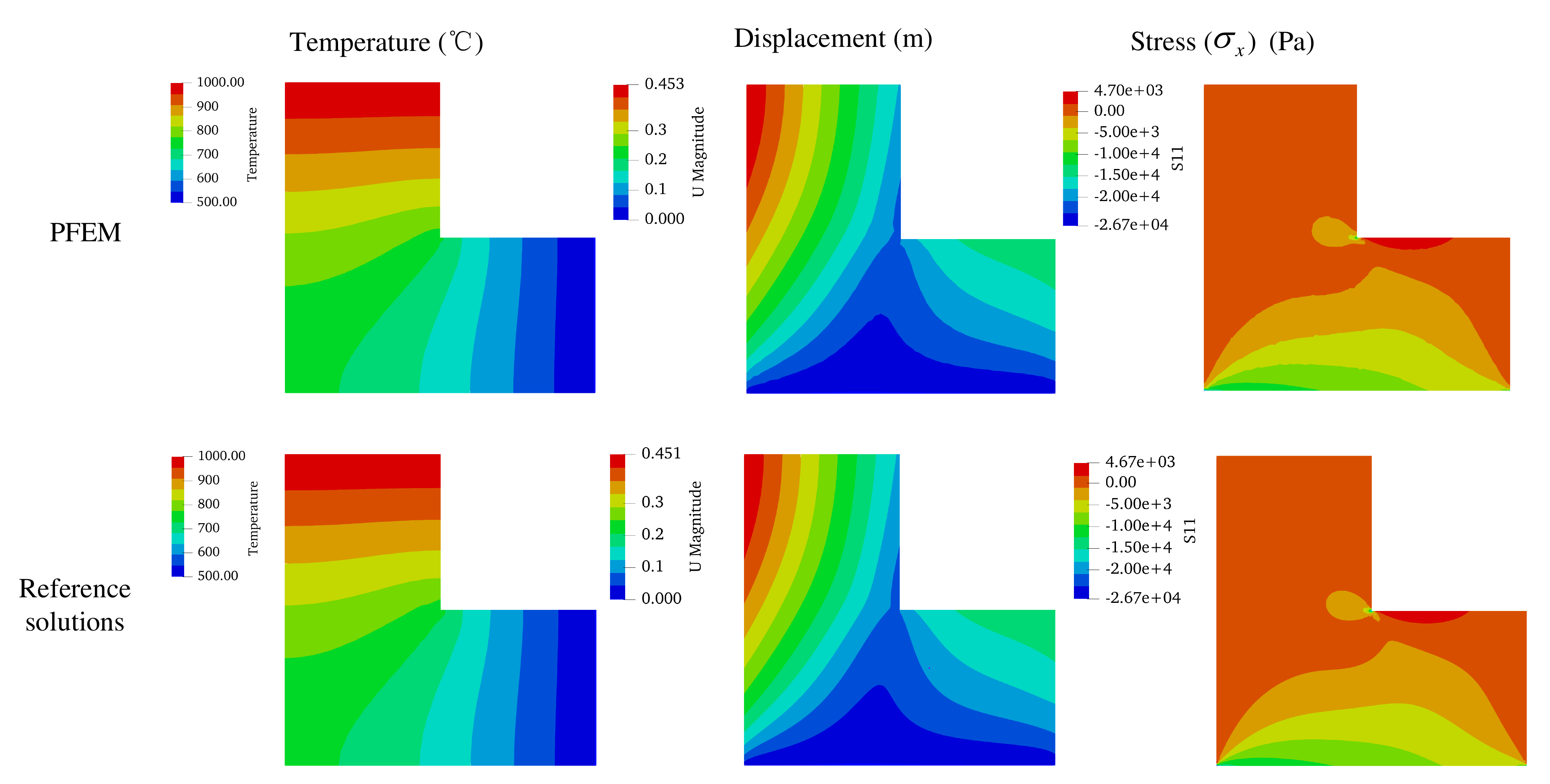}
  \caption{Contours for temperature, displacement and stress fields.}
  \label{fig:ex03_contour}
\end{figure}

As illustrated in Fig.~\ref{fig:ex03_error}, both the PFEM and conventional FEM exhibit satisfactory convergence behavior as the mesh is refined. Notably, the PFEM consistently achieves higher accuracy than FEM at equivalent mesh densities. The temperature and displacement histories at a designated monitoring point are shown in Figs.~\ref{fig:ex03_temp} and \ref{fig:ex03_disp}, respectively. Results from both the PFEM and FEM demonstrate excellent agreement with the reference solution. Furthermore, Fig.~\ref{fig:ex03_contour} presents the distributions of temperature, displacement, and horizontal stress at the final simulation time, further confirming the accuracy of the proposed PFEM in solving transient thermo-mechanical problems.

\subsubsection{Square panel with a complex cavity}
\begin{figure}[H]
  \centering
  \includegraphics[width=0.5\textwidth]{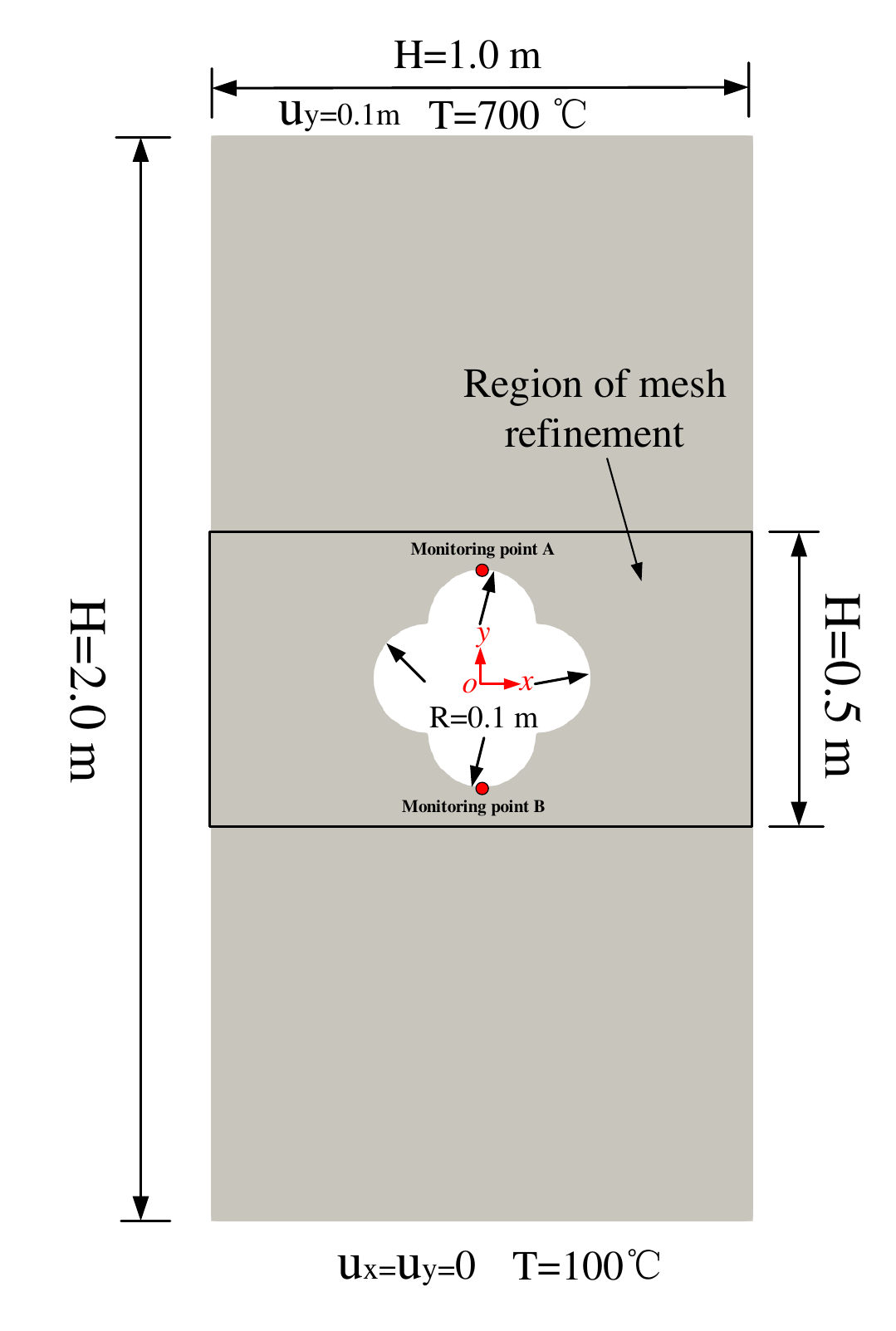}
  \caption{Schematic diagram of a square panel with multiple holes.}
  \label{fig:ex04Geo}
\end{figure}

To demonstrate the capability of the proposed PFEM in solving transient thermal stress problems, a rectangular body with a complex cavity is considered, as illustrated in Fig.~\ref{fig:ex04Geo}. The domain has a height of \( H = 2.0 \, \text{m} \) and a length of \( L = 1.0 \, \text{m} \), with the cavity formed by the intersection of four circles, each having a radius of \( 0.1 \, \text{m} \). The bottom surface is maintained at \( 100^\circ \text{C} \), while the top surface is heated to \( 700^\circ \text{C} \). A vertical displacement of \( 0.1 \, \text{m} \) is applied to the top surface. The material properties are as follows: thermal conductivity \( k = 80 \, \text{W/m·}^\circ\text{C} \), density \( \rho = 3870 \, \text{kg/m}^3 \), Young's modulus \( E = 1.2 \times 10^8 \, \text{Pa} \), Poisson’s ratio \( \nu = 0.3 \), thermal expansion coefficient \( \alpha = 0.0011 \, ^\circ \text{C}^{-1} \), and specific heat \( c = 0.45 \, \text{J/g·}^\circ\text{C} \). The total simulation time is set to 50~s with a uniform time step of 0.1~s.

\begin{figure}[H]
  \centering
  \includegraphics[width=0.9\textwidth]{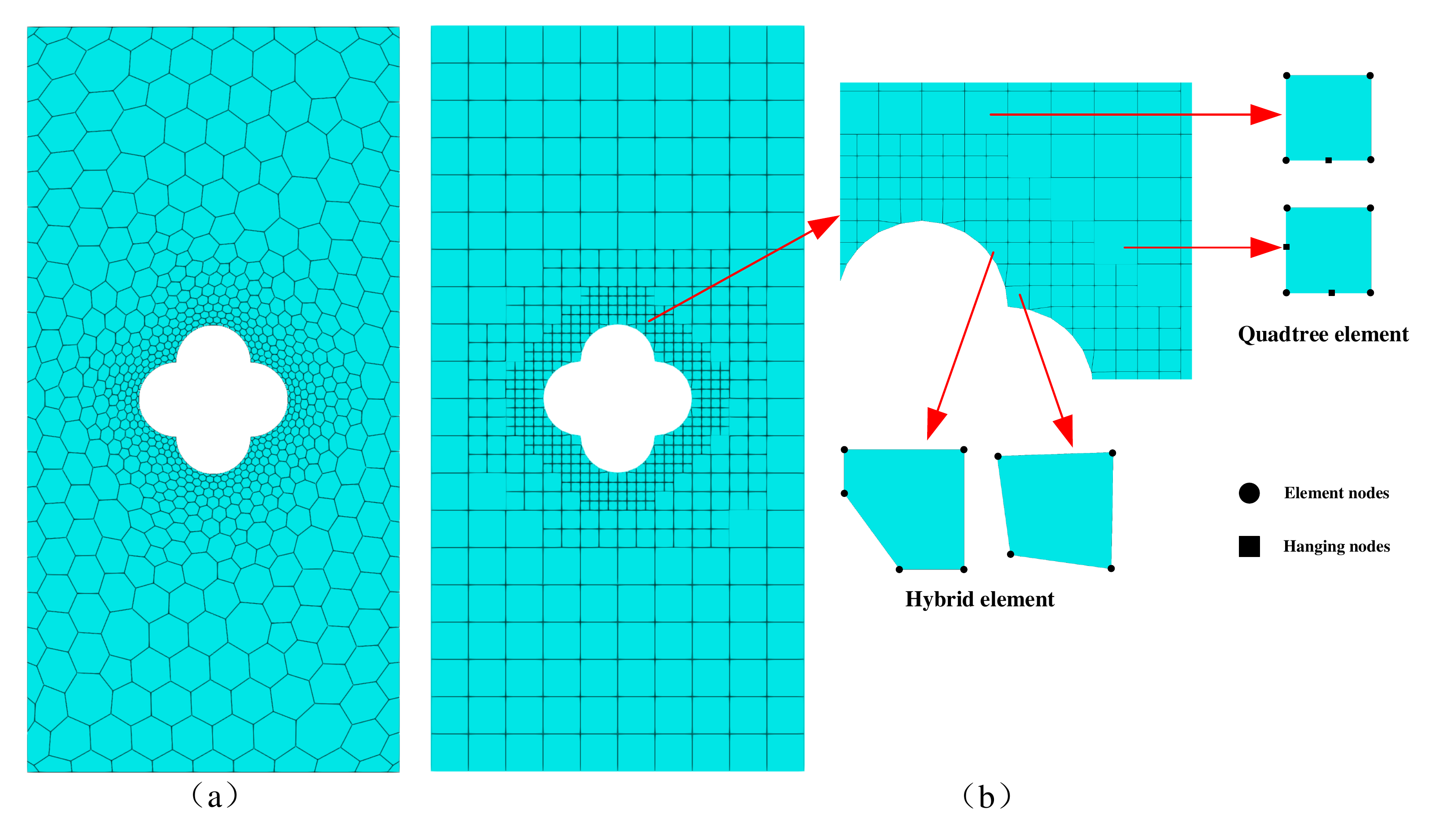}
  \caption{Meshes of square panel with multiple holes; (a) polygonal mesh (b) hybrid mesh.}
  \label{fig:ex04_mesh}
\end{figure}

The geometry is discretized using both polygonal and hybrid quadtree meshes with an element size of \( 0.1 \, \text{m} \), as shown in Fig.~\ref{fig:ex04_mesh}. In the hybrid quadtree mesh (Fig.~\ref{fig:ex04_mesh} b), regular quadtree elements are employed in simple regions, while irregular polygonal elements are used to conform to the complex cavity boundaries. The use of irregular elements enables the hybrid mesh to better capture geometric details along the boundary.

\begin{figure}[H]
  \centering
  \includegraphics[width=0.9\textwidth]{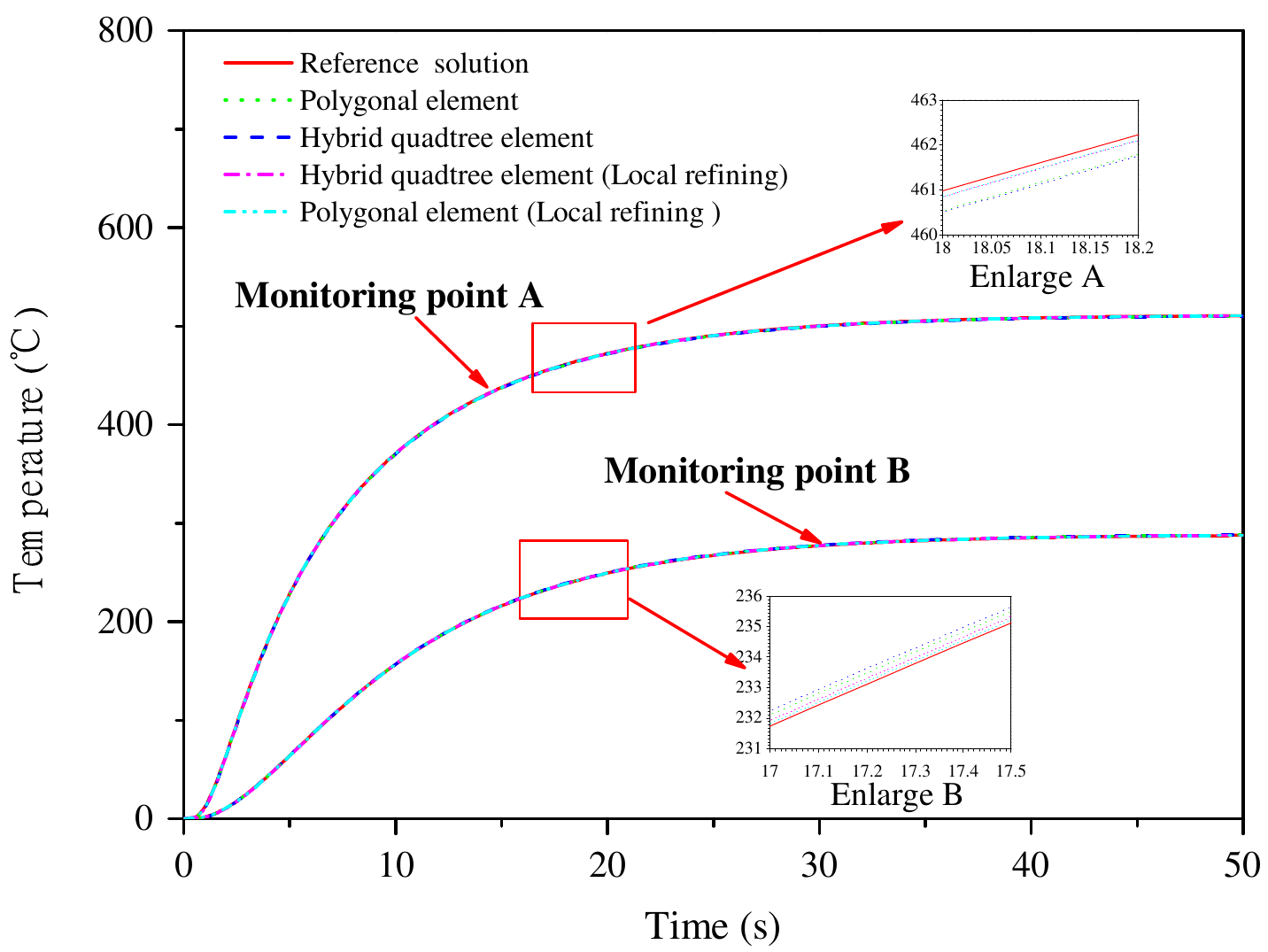}
  \caption{Comparison of the polygonal element and hybrid quadtree element in the history of temperature.}
  \label{fig:ex04_Temperature}
\end{figure}

\begin{figure}[H]
  \centering
  \includegraphics[width=0.9\textwidth]{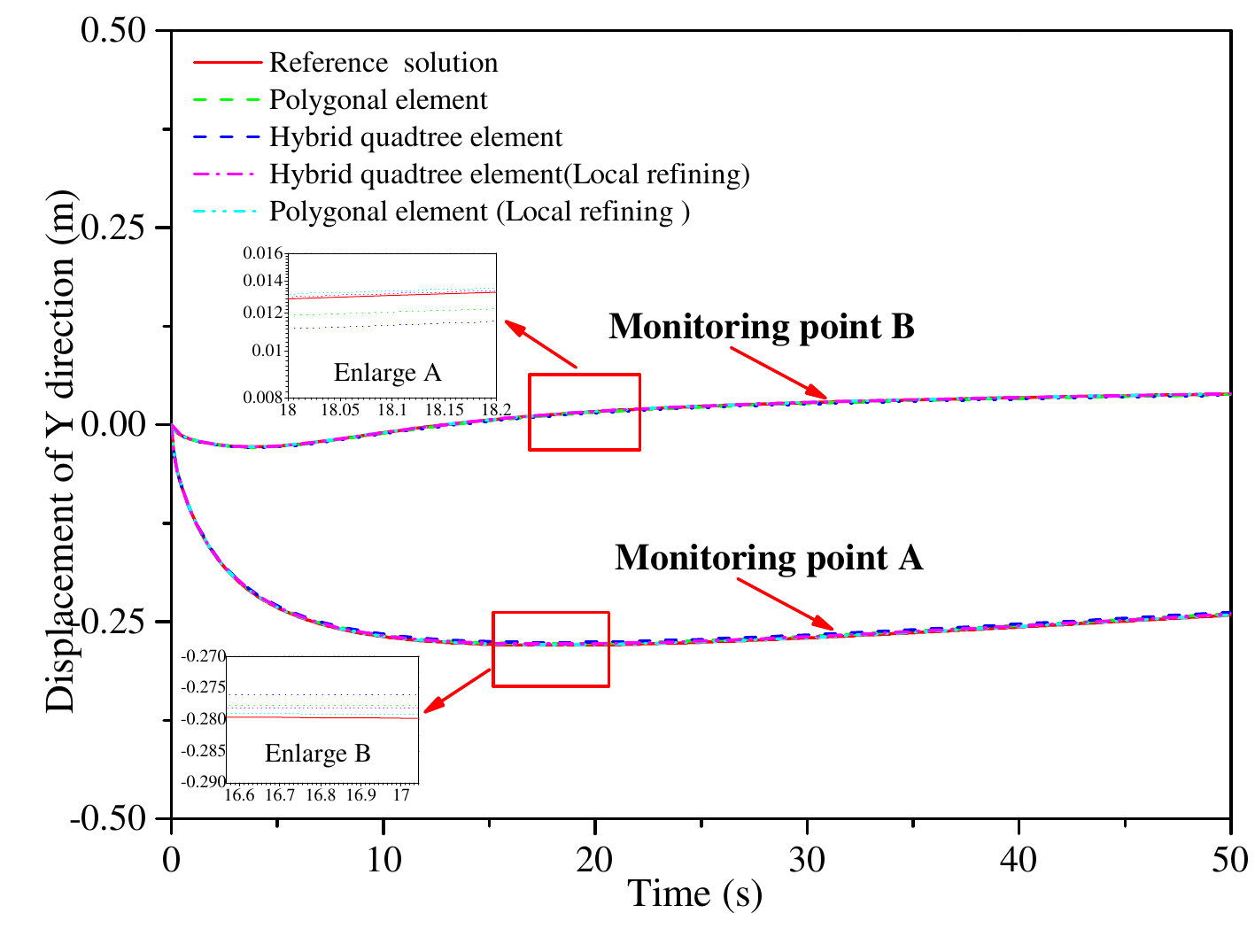}
  \caption{Comparison of the polygonal element and hybrid quadtree element in the history of displacement.}
  \label{fig:ex04_Displacement}
\end{figure}

\begin{figure}[H]
  \centering
  \includegraphics[width=1.0\textwidth]{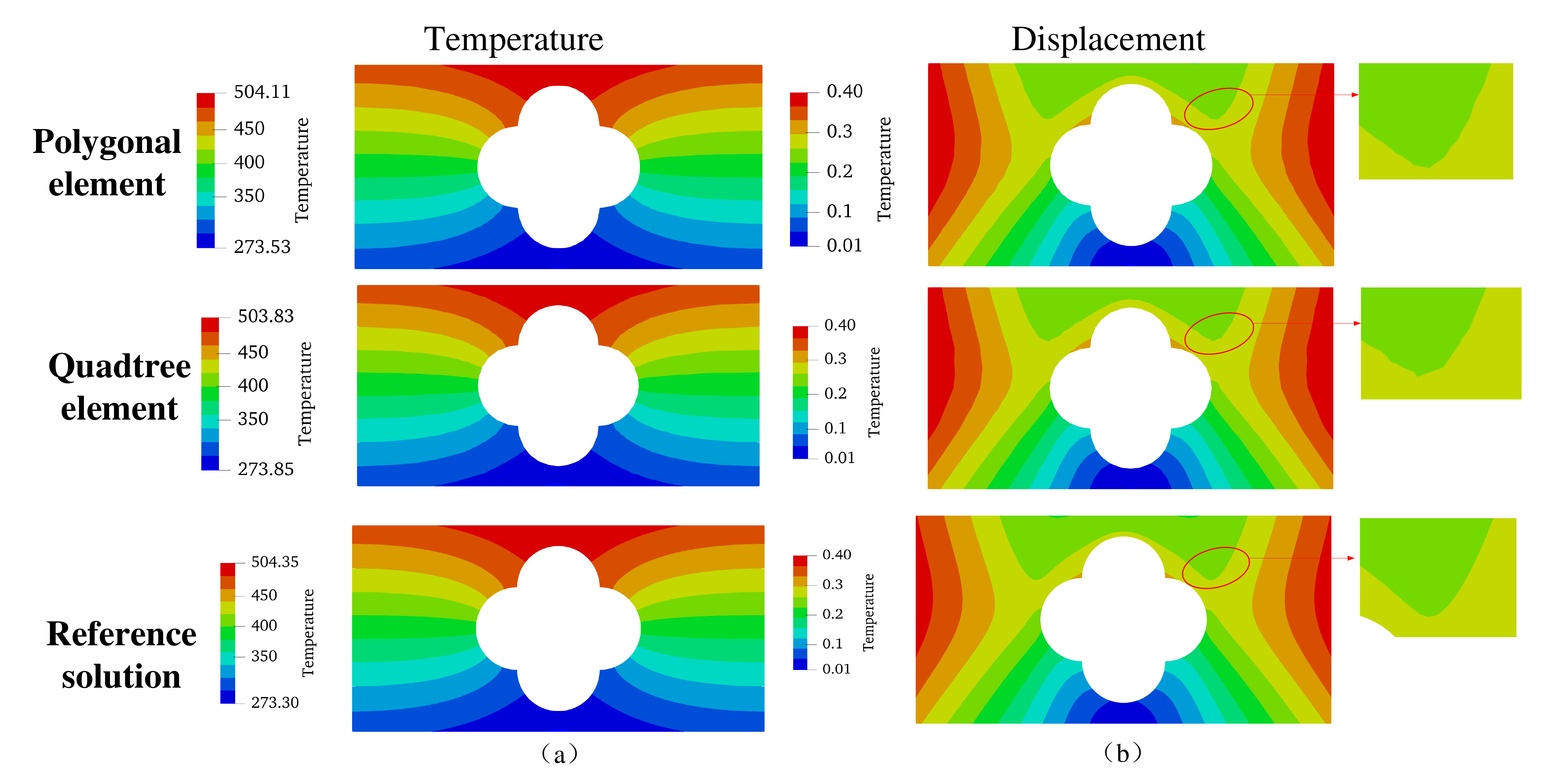}
  \caption{Contours before mesh refinement; (a) temperature distributions (Unit:$^\circ\mathrm{C}$); (b) displacement distributions (Unit:m).}
  \label{fig:ex04_contour}
\end{figure}

The temperature and displacement histories at monitoring points A and B are presented in Figs.~\ref{fig:ex04_Temperature} and \ref{fig:ex04_Displacement}, respectively. The results show good agreement with the reference solution. Fig.~\ref{fig:ex04_contour} illustrates the distributions of temperature and displacement. The overall distributions obtained using both polygonal and hybrid meshes are consistent with the reference. However, as seen in Fig.~\ref{fig:ex04_contour} (b), the displacement field obtained with hybrid quadtree elements exhibits less smoothness compared to that of the polygonal mesh. Quantitative comparisons in Tab.~\ref{tab:ex04_t1} show that polygonal elements yield lower relative errors than hybrid quadtree elements, which is attributed to the larger number of nodes per element in the polygonal mesh.

\begin{figure}[H]
  \centering
  \includegraphics[width=1.0\textwidth]{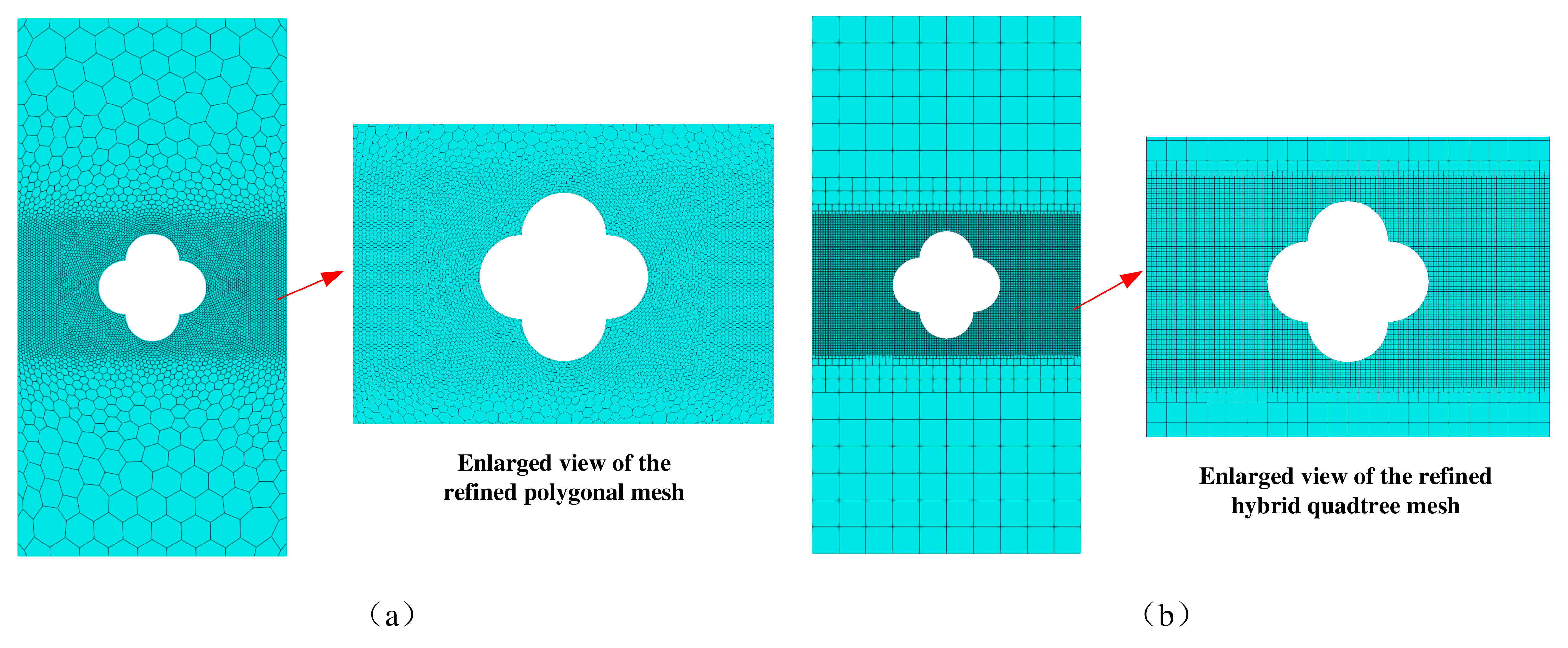}
  \caption{Local refining meshes; (a) polygonal mesh; (b) hybrid quadtree mesh.}
  \label{fig:ex04_mesh_refining}
\end{figure}

\begin{figure}[H]
  \centering
  \includegraphics[width=1.0\textwidth]{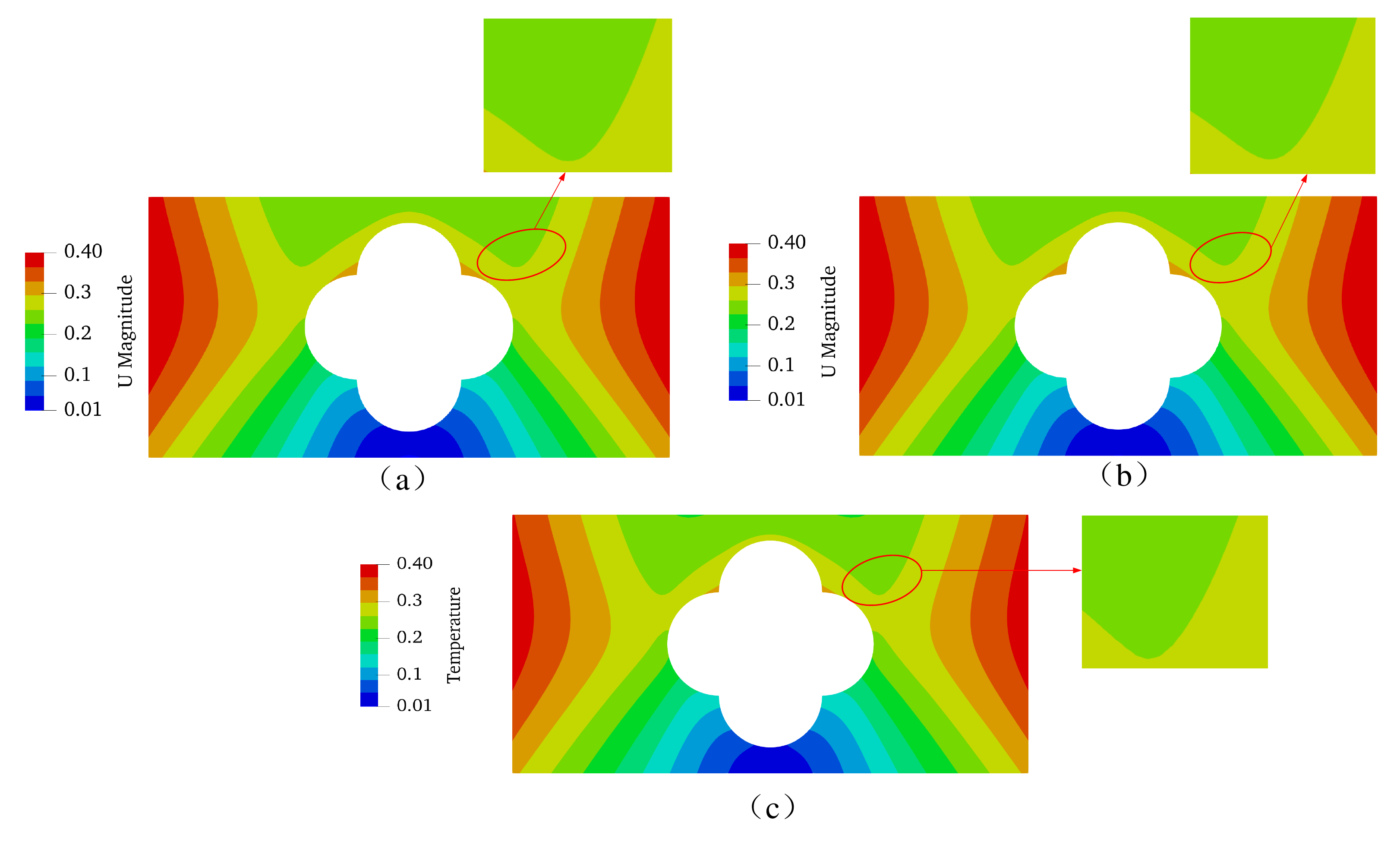}
  \caption{Displacement contour after mesh refinement; (a) polygonal elements; (b) hybrid quadtree mesh; (c) reference solution; (Unit:m)}
  \label{fig:ex04_contour_refining}
\end{figure}

To further improve the accuracy at the monitoring points, local mesh refinement is performed around those regions, as illustrated in Fig.~\ref{fig:ex04Geo}, and the refined mesh is shown in Fig.~\ref{fig:ex04_mesh_refining}. As indicated in Tab.~\ref{tab:ex04_t1}, local refinement significantly reduces the computational error at the monitoring points. The locally magnified views in Figs.~\ref{fig:ex04_Temperature} and \ref{fig:ex04_Displacement} also confirm the improved accuracy due to mesh refinement. Moreover, Fig.~\ref{fig:ex04_contour_refining} demonstrates that the displacement and stress contours become notably smoother after refinement, indicating enhanced solution quality.

Notably, after refinement, the number of nodes in the polygonal and hybrid quadtree meshes becomes comparable, while the polygonal mesh uses 49.5\% fewer elements. Consequently, the computational time for the refined polygonal mesh is shorter than that of the hybrid quadtree mesh. However, Tab. \ref{tab:ex04_t1} also show that the computational time is significantly reduced through acceleration technology for the quadtree mesh.

\begin{table}[H]
\centering
\caption{Relative errors for the monitoring points.}
\resizebox{\textwidth}{!}{ 
\begin{tabular}{ccccccc}
\toprule
Element type  & Elements & Nodes & Relative error for A &  Relative error for B & CPU time (s)  & CPU time (s)*\\ 
\midrule
Polygonal element            & 707       & 1422        & 8.21$\times10^{-4}$  & 1.48$\times10^{-3}$  & 179.30 & —\\ 
Hybrid quadtree element      & 536       & 630      & 9.78$\times10^{-4}$   & 1.87$\times10^{-3}$  & 112.20 & 63.17\\ 
Refining polygonal element   & 5715     & 11438     & 2.49$\times10^{-4}$   & 3.68$\times10^{-4}$ &  1614.30 & —\\ 
Refining hybrid quadtree element     & 11328      & 11664  & 3.09$\times10^{-4}$  & 6.99$\times10^{-4}$ &  2181.60 & 1346.05\\ 
\bottomrule 
\label{tab:ex04_t1}
\end{tabular}}
\vspace{-2em}
\begin{flushleft}
\small Note: * denotes the use of acceleration technology based on the parent element.
\end{flushleft}
\end{table}

\section{Conclusions}
\label{sec:conclusions}
This paper proposed a PFEM for analyzing steady-state and transient thermal stress problems in two-dimensional continua. Wachspress rational basis functions were adopted to construct conforming interpolation over convex polygonal meshes, enabling flexible geometric modeling and improved numerical accuracy. The main findings are summarized as follows:

(1) The PFEM formulation naturally supports non-matching meshes and hanging nodes without requiring additional constraints, making it well-suited for adaptive meshing and multi-scale applications.

(2) Compared to conventional FEM, the PFEM demonstrates superior convergence behavior and solution accuracy, particularly in scenarios involving local mesh refinement or geometric irregularities.

(3) The quadtree-based acceleration technique, which reuses precomputed stiffness and mass matrices, significantly reduces computational cost while maintaining high solution fidelity.

Future research could focus on enhancing the performance of PFEM for higher-dimensional and nonlinear thermal conduction problems, as well as exploring its applications in real-world engineering contexts.

\section{Acknowledgments}
The Xing Dian Talent Support Program of Yunnan Province (grant NO. XDYC-QNRC-2022-0764), the Yunan Funndamental Research Projects (grant NO. 202401CF070043) and the Science and Technology Talents and Platform Plan (grant NO. 202305AK34003) provided support for this study.
\appendix

 \bibliographystyle{elsarticle-num} 
 \bibliography{cas-refs}





\end{document}